\documentclass[11pt]{amsart}
\usepackage{amscd,amssymb}
\usepackage{amsthm,amsmath,amssymb}
\usepackage{youngtab}
\usepackage{tikz}
\usepackage{graphicx}
\usepackage{color}
\usepackage{pdflscape}
\usepackage{bm}
\usepackage{mathrsfs}
\usepackage{verbatim}
\usepackage{url}
\usepackage{xtab}
\usepackage{longtable}
\usepackage{enumitem}    

\setlist[itemize]{leftmargin=22pt}   

\newtheorem{lemma}{Lemma}[section]
\newtheorem{theorem}[lemma]{Theorem}
\newtheorem*{theorem*}{Theorem}
\newtheorem*{conjecture*}{Conjecture}

\newtheorem*{proposition*}{Proposition}
\newtheorem{proposition}[lemma]{Proposition}
\newtheorem{corollary}[lemma]{Corollary}

\theoremstyle{definition}
\newtheorem{definition}[lemma]{Definition}
\newtheorem*{definition*}{Definition}

\newtheorem{example}[lemma]{Example}

\newtheorem{remark}[lemma]{Remark}

\newtheorem{acknowledgements*}[lemma]{Acknowledgements}

\makeatletter
\providecommand{\leftsquigarrow}{%
  \mathrel{\mathpalette\reflect@squig\relax}%
}
\newcommand{\reflect@squig}[2]{%
  \reflectbox{$\m@th#1\rightsquigarrow$}%
}
\makeatother

\newcommand{\ra}{\rightarrow}

\newcommand{\llra}{\longleftrightarrow}

\newcommand{\ob}{\overline}
\newcommand{\Ra}{\Rightarrow}

\newcommand{\ep}{\varepsilon}

\newcommand{\al}{\alpha}
\newcommand{\be}{\beta}
\newcommand{\ga}{\gamma}
\newcommand{\Ga}{\Gamma}
\newcommand{\si}{\sigma}
\newcommand{\lam}{\lambda}

\newcommand{\ze}{\zeta}
\newcommand{\de}{\delta}
\newcommand{\De}{\Delta}

\newcommand{\Up}{\Upsilon}
\newcommand{\IH}{\mathcal{H}}

\newcommand{\Z}{\mathbb{Z}}

\newcommand{\Area}{\textup{Area}}
\newcommand{\rArea}{\textup{A}\widetilde{\textup{re}}\textup{a}}

\newcommand{\ms}{\medskip}

\newcommand{\noi}{\noindent}
\newcommand{\qq}{\qquad}
\newcommand{\qqq}{\qq \qq \qq \qq}
\newcommand{\sm}{\!\!\!\!}
\newcommand{\les}{\leqslant}
\newcommand{\ges}{\geqslant}

\newcommand{\eq}{\,=\,}
\newcommand{\qs}{\scalebox{0.5}{\,}}    
\newcommand{\qc}{\scalebox{0.25}{\!}}  
\newcommand{\pr}{\operatorname{pr}}
\newcommand{\rep}{\operatorname{rep}}
\newcommand{\orb}{\operatorname{Orb}}
\newcommand{\sh}{\operatorname{sh}}

\newcommand{\sgn}{\operatorname{sgn}}

\newcommand{\rdes}{\mathcal{R}}
\newcommand{\wei}{\mathscr{L}}
\newcommand{\rxi}{\ob{\rdes}^{\raise - 2.0pt \hbox{\scriptsize $\wei$}}}
\newcommand{\Sxi}{\ob{S}^{\raise - 2.0pt \hbox{\scriptsize $\wei$}}}
\newcommand{\nchooseq}{n \choose {\raise 1.6pt \hbox{\scriptsize $q$}} }
\newcommand{\stackclose}[2]{\stackrel{\scriptsize \raise -0.5pt \hbox{$#1$}}{#2}}
\newcommand{\app}{\approx}
\newcommand{\napp}{\not\approx}

\newcommand{\tskip}[1]{ \hskip - #1 pt \begin{tabular}{r p{{325 - #1} pt}} }
\newcommand{\udots}{ \raise 5pt \hbox{\reflectbox{$\ddots$}}}

\newcommand{\sdf}{\scalebox{1.5}}
\newcommand{\sdg}{\scalebox{1}}

\newcommand{\diagramB}{\begin{center}
\begin{picture}(180,28)
\put(  -3, 10){$B_n$}
\put(  -1, -2.7){$\scriptstyle{\wei}$}
\put( 30, 10){\circle*{5}}
\put( 30, 10){\line(1,0){30}}
\put( 42.5, 13){$\scriptstyle{4}$}
\put( 60, 10){\circle*{5}}
\put( 60, 10){\line(1,0){30}}
\put( 90, 10){\circle*{5}}
\put( 90, 10){\line(1,0){20}}
\put(120, 10){\circle*{1}}
\put(130, 10){\circle*{1}}
\put(140, 10){\circle*{1}}
\put(150, 10){\line(1,0){20}}
\put(170, 10){\circle*{5}}
\put( 28, 18.5){$t$}
\put( 57, 19){$s_1$}
\put( 87, 19){$s_2$}
\put(167, 19){$s_{n-1}$}
\put(28,-3){$\scriptstyle{b}$}
\put(58,-3){$\scriptstyle{a}$}
\put(88,-3){$\scriptstyle{a}$}
\put(168,-3){$\scriptstyle{a}$}
\end{picture}
\end{center}}

\newcommand{\arr}[2]{
\draw[->] (#1,#2)  to  [out=-70, in=70] (#1,#2-0.45);
}

\newcommand{\art}[2]{
\draw[->] (#1,#2)  to  [out=110, in=245] (#1,#2+2.3);
}

\newcommand{\argh}[2]{
\arr{#1}{#2}
\arr{#1}{#2-1}
\arr{#1+0.03}{#2-1.88}
\art{#1-0.25}{#2-2.4}
}

\newcommand{\arrk}[2]{
\draw[->] (#1-0.05,#2+0.35)  to  [out=-70, in=70] (#1-0.05,#2-0.42);
}

\newcommand{\artk}[2]{
\draw[->] (#1,#2)  to  [out=105, in=250] (#1,#2+2.7);
}

\newcommand{\arghk}[2]{
\arrk{#1}{#2}
\arrk{#1}{#2-1}
\arrk{#1}{#2-2}
\artk{#1-0.2}{#2-2.38}
}

\newcommand{\tbox}[4]{
\draw (#1,#2) -- (#1+#3,#2) -- (#1+#3,#2-#3) -- (#1,#2-#3) -- (#1,#2);
\node at (#1+0.5*#3,#2-0.5*#3) {$#4$};
}
\newcommand{\xbox}[4]{
\draw (#1,#2) -- (#1+#3,#2) -- (#1+#3,#2-#3) -- (#1,#2-#3) -- (#1,#2);
\node at (#1+0.54*#3,#2-0.54*#3) {$#4$};
}
\newcommand{\zbox}[4]{
\draw (#1,#2) -- (#1+#3,#2) -- (#1+#3,#2-#3) -- (#1,#2-#3) -- (#1,#2);
\node at (#1+0.5*#3,#2-0.4*#3) {$#4$};
}

\address{Department of Mathematics \\ National University of Singapore \\ 10 Lower Kent Ridge Road \\ Singapore 119076}
\email{edmund.howse@nus.edu.sg}
\title{Vogan classes in type $B_n$}
\author{Edmund Howse} 
\thanks{The author is supported by Singapore MOE Tier 2 AcRF MOE2015-T2-2-003. }
\begin{document}



{\color{white} $\bullet$}\\[-3.55em]

\begin{abstract}
\qc \qc  Kazhdan and Lusztig have shown how to partition a Coxeter group into cells. In this paper, we use the theory of Vogan classes to obtain a first characterisation of the left cells of type $B_n$ with respect to a certain choice of weight function. 
\end{abstract}

\maketitle

\section{Introduction}\label{intro}


Lusztig has described how to partition a Coxeter group into left, right and two-sided cells with respect to a weight function \cite{LWeyl}. This is done via certain equivalence relations that are calculated in the corresponding Iwahori--Hecke algebra, and the resulting cells afford representations of both the group and the algebra. Algebraic techniques have been developed to reduce the determination of cells to combinatorial calculations at the level of the group. 

These ideas, along with their connections to other areas of mathematics, were first outlined by Kazhdan and Lusztig \cite{KL}. Their paper contains a quintessential result in the theory of cells -- the classification of cells in type $A_n$. In this setting,
the left cells afford a complete list of irreducible representations of the corresponding Iwahori--Hecke algebra, a pair of left cells afford isomorphic representations if and only if they are contained in the same two-sided cell, and two elements of the group are in the same left cell if and only if they have the same recording tableaux under the Robinson--Schensted correspondence.

There is ample motivation then, to have an interest in the theory of cells. One of the aims of this theory is to classify the cells of finite Coxeter groups.

Only three types of finite irreducible Coxeter groups may have associated Iwahori--Hecke algebras with unequal parameters; these are type $F_4$, type $B_n$, and type $I_2(m)$, where $m$ is even. 
The theory here is more complex than the equal parameter case; the conditions for determining whether two weight functions are cell-equivalent -- that is, they give rise to the same partition of the Coxeter group into cells -- are far from obvious.

Lusztig comprehensively discussed the cells of type $I_2(m)$ in \cite{L}, while Geck used a combination of theoretical considerations and explicit computer calculations to resolve the case of type $F_4$ in \cite{GeckF}. This leaves the case of type $B_n$ to be considered.

Denote by $W_n$ the Coxeter group of type $B_n$. The weighted Coxeter system $(W_n,S,\wei)$ may be concisely described by its Coxeter diagram:
\diagramB
{\color{white} $\bullet$}\\[-1.8em]

When $b/a = 1$, the cells of $W_n$ are known following a series of papers of Garfinkle \cite{GarI}, \cite{GarII}, \cite{GarIII}. In this case, the left cells may be combinatorially described using domino tableaux and associated operations.
 A description of the left cells when $b/a \in \{\frac12,\frac32\}$ is due to Lusztig \cite{LWeyl}. 
 The other known case is detailed in papers of Bonnaf\'{e} and Iancu \cite{BI} and Bonnaf\'{e} \cite{BonAsy}; any weight function such that $b/a > n-1$ corresponds to an ``asymptotic'' choice of parameters, and the resulting cells are governed by a generalised Robinson--Schensted correspondence. The aforementioned results of Kazhdan and Lusztig on cells in type $A_n$ all have an analogue in type $B_n$ in the asymptotic case, but not in the equal parameter case. 

An important development in the study of the cells of $W_n$ came in the form of a number of conjectures by Bonnaf\'e, Geck, Iancu and Lam \cite{BGIL}. These conjectures state conditions for two weight functions on $W_n$ to be cell-equivalent, as well as a unified combinatorial description of the left, right and two-sided cells for each of these cases. Although there are results in this direction due to Bonnaf\'e \cite{BonKLB} \cite{BonKLBErratum}, a proof of these conjectures remains elusive. 

In this paper, our focus is on the left cells of $W_n$ 
 when $b/a > n-2$. The conjectures of \cite{BGIL} suggests that the corresponding weight functions belong to one of the following three classes of cell-equivalence:

\begin{itemize}
\item when $b/a > n-1$; this is the {\em asymptotic case},
\item when $b/a = n-1$; which we refer to as the {\em intermediate case},
\item when $b/a \in (n-2,n-1)$; referred to as the {\em sub-asymptotic case}.
\end{itemize}

In the following two sections we recall some necessary background material. In Section 4 we generalise the notion of the {\em enhanced right descent set} of Bonnaf\'e and Geck \cite{BG} to obtain an invariant $\rxi$ of left cells for finite weighted Coxeter systems. Section 5 is a collection of technical results for later use. In Section 6, we recall the concept of Vogan classes from \cite{BG}, and establish a particular set $\Xi$ of KL-admissible pairs for use in type $B_n$, valid exactly when $b/a > n-2$. 
These KL-admissible pairs describe maps that can be used to determine 
cellular information that is common to the three cases mentioned above; Section 7 is dedicated to understanding this.
In the final section,  we determine the left Vogan $(\Xi,\rxi)$-classes, which leads to the following result. 

\begin{theorem*}
Suppose $(W_n, S, \wei)$ is such that $b/a \ges n-1$. Then two elements of $W_n$ 
 are in the same left cell if and only if they lie in the same left Vogan $(\Xi,\rxi)$-class.
\end{theorem*}

As such, the left Vogan $(\Xi,\rxi)$-classes offer a new characterisation of the asymptotic left cells of $W_n$, 
 and the first characterisation of the left cells when $b/a = n-1$. 


\section{Kazhdan--Lusztig cells with unequal parameters}

Let $(W,S)$ be a Coxeter system, let $y, w \in W$, $s,t \in S$, and denote by $\ell:W\ra\Z_{\ges0}$ the standard length function. A {\em weight function} for $W$ is any map $\wei:W\ra \Z$ such that $\wei(yw) = \wei(y) + \wei(w)$ whenever $\ell(yw) = \ell(y) + \ell(w)$. A weight function is uniquely determined by its values on $S$; 
conversely, any function $\wei:S\ra\Z$ such that $\wei(s)=\wei(t)$ whenever 
$st$ has odd order extends uniquely to a weight function on $W$. 

We will assume throughout that $\wei(s)>0$ for all $s \in S$; see Lusztig \cite{LWeyl} and Bonnaf\'{e} \cite{BonSP} for the original, more general framework. 

\ms

We will denote by $\les$ the Bruhat--Chevalley order on $W$, and write $y < w$ if $y \les w$ with $y \neq w$. If $y$ is a suffix of $w$, we write $y \les_e w$; see 2.1 of Geck--Pfeiffer \cite{GP} for details. 

\ms

Let $I \subseteq S$ be non-empty. The parabolic subgroup $W_I:=\langle I \rangle$ has a corresponding set $X_I$ of distinguished left coset representatives. 
For all $w \in W$, there exist unique elements $x \in X_I$, $u \in W_I$ such that $w=xu$; moreover, $\ell(w) = \ell(x) + \ell(u)$. In this context, we denote:
\[\rep_I(w):=x,\ \ \ \pr_I(w):=u.\]
There exists a bijection:
\[ \begin{array}{ccc} W & \llra & X_I\times W_I, \\ w & \llra & (\rep_I(w),\, \pr_I(w)). \end{array}\]

\smallskip

The Iwahori--Hecke algebra $\IH:=\IH(W,S,\wei)$ is a deformation of the group algebra of $W$ over $\mathcal{A}:=\Z[v,v^{-1}]$, the ring of Laurant polynomials with indeterminate $v$. 

The Iwahori--Hecke algebra has an $\mathcal{A}$-basis $\{T_w \,:\, w\in W\}$; multiplication between basis elements may be described by the formula:
\[ T_sT_w = \left\{ \begin{array}{ll} T_{sw} & \ \  \textrm{ if } sw>w, \\  T_{sw} + (v^{\wei(s)} - v^{-\wei(s)})\,T_w  & \ \ \textrm{ if } sw<w. \end{array} \right.  \]

We refer to the elements of the set 
$\{v^{\wei(s)} \, : \, s \in S\}$ as the {\em parameters} of the Iwahori--Hecke algebra. If there exist $s,t\in S$ such that $v^{\wei(s)} \neq v^{\wei(t)}$ then we say that $\IH$ has {\em unequal parameters}. Otherwise, we say that we are in the {\em equal parameter} case. 

\ms

Multiplication in the Kazhdan--Lusztig basis $\{ C_w \,: \, w \in W\}$ of $\IH$ (see Kazhdan--Lusztig \cite{KL} and Lusztig \cite{LWeyl}) determines the left cells of $W$ as follows. The relation defined by:
\[
y\les'_{L}w \ \text{ if } \ \left\{\!\begin{array}{ll} \text{there exists } s \in S \text { such that }  \\ C_y \text{ occurs with non-zero coefficient in } C_sC_w \end{array}\right.
\]


\noi can be extended to its reflexive, transitive closure -- a preorder $\les_L$ called the {\em Kazhdan--Lusztig preorder}. The associated equivalence relation on $W$, denoted $\sim_L$, is defined by:
\[y \sim_L w \iff  y \les_L w \,\text{ and }\, w \les_L y.\] 
The resulting equivalence classes are called {\em left cells}. 
Similar definitions exist for the preorder $\les_R$ and right cells, as well as the preorder $\les_{LR}$ and two-sided cells; see Lusztig \cite{L} \cite{Lupdated} for details. 


The Coxeter group $W_I$ has its own left, right and two-sided cells, arising from relations denoted $\sim_{L,I}$, $\sim_{R,I}$ and $\sim_{LR,I}$ respectively.

\section{The Coxeter group of type $B_n$}\label{B_n section}

We retain the notation for the Coxeter group of type $B_n$ from \S\,\ref{intro}.
It is often useful to identify $W_n$ with the group of signed permutations. The map given by $ t \mapsto (1,-1)$, $s_i \mapsto (i,i+1)(-i,-i-1)$
defines an isomorphism between these groups.  Thus we can write $w \in W_n$ as a sequence $w(1), \ldots, w(n)$, where for $1 \les i \les n$ we have $w(i) = \ep_i p_i$, with $\ep_i \in \{\pm 1\}$, and $p_1, \ldots, p_n$ forming a permutation of $n$.

 One benefit of this identification is the use of the following classical result. Let $1 \les i \les n-1$ and $1 \les j \les n$. Set $t_j:=s_{j-1} \cdots s_1 t s_1 \cdots s_{j-1} = (j,-j)$.

\begin{lemma}\label{signdes} Let $w \in W_n$. Then:
\smallskip

\setlength{\tabcolsep}{4pt}  
\hskip -5.5pt
\begin{tabular}{rl}
\sm (i) &  $\ell(ws_i) < \ell(w) \iff w(i+1) < w(i)$, \\[0.1em]
(ii) & $\ell(wt_j) < \ell(w) \iff w(j)<0$.
\end{tabular}

\end{lemma}


There exists a generalised Robinson--Schensted correspondence from the elements of $W_n$ to pairs of standard bitableaux of size $n$ with the same shape; see Bonnaf\'e--Iancu \cite{BI} for details. We denote this correspondence
\[ w\, \longmapsto \big( A_n(w),\, B_n(w)\big),\]

\noi and set $\sh(w)$ to be the shape of $A_n(w)$ which is equal to the shape of $B_n(w)$. The bitableau $A_n(w)$ is called the {\em insertion bitableau of $w$}, and $B_n(w)$ is called the {\em recording bitableau of $w$}.

\ms

In \S\,1 of \cite{BI}, it is noted that the left cells of $(W_n,S,\wei)$ are independent of the exact value of $b/a$, provided it is sufficiently large (with respect to $n$). Such a weight function is termed an {\em asymptotic weight function}. Whenever $W_n$ is equipped with an asymptotic weight function, we say that we are in the {\em asymptotic case}.

\begin{theorem}\label{ThmBI}
(Bonnaf\'e--Iancu,  Theorem 7.7 of \cite{BI}; Bonnaf\'{e}, Remark 3.7  and Proposition 5.5 of \cite{BonAsy})  A weight function on $W_n$ is asymptotic if and only if $b/a > n-1$. Moreover, let $y,w \in W_n$ and suppose that we are in the asymptotic case. Then:  

\smallskip

\setlength{\tabcolsep}{4pt} \hskip -5.5pt \begin{tabular}{rl}
(i) &  $y \sim_L w\iff B_n(y) =  B_n(w)$, \\[0.1em]
(ii) & $y \sim_R w\iff A_n(y) = A_n(w)$, \\[0.1em]
(iii) & $y \sim_{LR} w \iff \sh(y) = \sh(w)$. 
\end{tabular}

\end{theorem}

We conclude with some additional notation. Denote by $\ell_t$ the function that counts the number of occurences of the generator $t$ in a reduced expression for $w \in W_n$. Let $e$ be the identity element, and let $w_0$ be the longest word of $W_n$.


\section{A new descent set for finite Coxeter groups}
Let $(W,S,\wei)$ be a weighted Coxeter system. The {\em right descent set of $w \in W$} is given by
\[\rdes(w):=\{\, s \in S \,:\, \ell(ws)<\ell(w)\,\}.\]
The right descent set is an invariant of the left cells of $W$ (see \cite{LWeyl}), so if $y\sim_L w$ then $\rdes(y)=\rdes(w)$. This concept may be refined as in Example 7.4 of Bonnaf\'{e}--Geck \cite{BG}. Let
 \[S^\wei:=\,S \,\cup\,\{\,sts \,:\, s,t \in S\text{, such that } \wei(t)>\wei(s)\,\}.\]
Then the {\em enhanced right descent set} of $w \in W$ is
\[\rdes^{\raise  0.5pt \hbox{\scriptsize $\wei$}}(w):=\qs \{ \,\si \in S^\wei \,:\, \ell(w \si)<\ell(w)\,\}.\]
This is again an invariant of left cells. Although this is only a slightly finer invariant than the right descent set, its strength lies in being fine enough to determine left cells of dihedral groups with respect to any of the three cell-equivalence classes of weight functions. It is therefore a useful tool when considering parabolic subgroups of $W$ of rank $2$. In this section, we further refine this concept to obtain an invariant of left cells of finite weighted Coxeter systems that is sensitive to the choice of weight function.

\subsection{A generalised enhanced right descent set}

\begin{definition}\label{SPD}
Let $(W,S,\wei)$ be a weighted Coxeter system, and let
\[ \Sxi:= S \cup \left\{ \,  s_{k} \cdots s_1 t s_1 \cdots s_{k} \ :\, \!\begin{array}{ll} \wei(t) > k\cdot \wei(s_i) \ \text{ and } \\ \textup{order}({s_is_{i+1}})=3 \text{ for } 1 \les i \les k-1 \end{array} \! \right\}\qc.\]

\noi For $w \in W$, let  \\[-1em]
\[\rxi(w):=  \{\, \si \in \Sxi \,:\, \ell(w \si) < \ell(w) \,\}.\]
\end{definition}

\noi Note that $\rxi(w) = \rdes(w)$ when $\wei$ is constant on $S$.

\begin{lemma}\label{F_4 prop}
Let $W$ be of type $F_4$, and let $\wei:W\ra \Z$ be 
 a weight function. Then $\rxi$ is an invariant of the left cells of $W$ with respect to $\wei$.
\end{lemma}

Proof -- The cell-equivalence classes of weight functions on $W$ are known following Corollary 4.8 of Geck \cite{GeckF}. The result can then be verified for a representative of each class with some elementary computer code in conjunction with the Python module PyCox; see Geck \cite{PyCox}. \qed

\begin{remark}\label{SPD W_n ver.}
For $W=W_n$ and $b/a \ges 1$, we have \[\Sxi \subseteq \qs \{s_1,\, s_2,\ldots, s_{n-1},\, t_1,\, t_2, \ldots, t_n\},\] and for $w \in W_n$, we have \\[-1em]
\[\rxi(w) \eq \rdes(w) \cup\{\, t_k \ : \  \frac{b}a > k-1  \ \textup{ and }\  \ell(wt_k) < \ell(w)\,\}.\]
We set $K:=\{t,s_1,\ldots,s_{n-2}\}$ so that $W_K = W_{n-1}$. We may describe the set $X_K$ as the set of all suffixes of the Coxeter word $t_n = s_{n-1} \cdots s_1ts_1 \cdots s_{n-1}$. 
\end{remark}

From now until the end of this section only, we adopt the setup of Remark~\ref{SPD W_n ver.}.

\begin{remark}\label{incl}
Let $x \in X_K$ and $u \in W_K$. Then $\rxi(xu) = \rxi(u)\, \cup\, A$, for some 
$A \subseteq \{s_{n-1},\, t_n\}$.
\end{remark}

\begin{proposition}\label{SPD invariant}
Let $y,w \in W_n$. If $y \sim_L w$, then $\rxi(y)=\rxi(w)$.
\end{proposition}

Proof -- Fix a weight function $\wei$ as in Remark~\ref{SPD W_n ver.}, let $y,w \in W_n$, and suppose that $y \sim_L w$. Recall that if $y \sim_L w$ then $\rdes(y)=\rdes(w)$, and so we only need to evaluate the membership of $t_j$  in $\rxi(y)$ and $\rxi(w)$, where $j$ satisfies both $2 \les j \les n$ and $b/a > j-1$. Proceed by induction on $n$.

\ms

Suppose $n=2$. Then $\rxi$ coincides with $\rdes^{\raise - 0.5pt \hbox{\scriptsize $\wei$}}$, and $W_2$ coincides with $W(I_2(4))$. The statement is then true by Example 7.4 of \cite{BG}. Now assume the statement is true for all $r < n$. We show that it is true for $r=n$.

\ms

Let $2 \les j \les n-1$ and suppose $b/a > j-1$. Then by Remark~\ref{incl}, we have
\begin{equation}\label{eqn 01}
t_j \in \rxi(w) \iff t_j \in \rxi(\pr_K(w)).
\end{equation}

By Theorem 1 of Geck \cite{GeckInd}, $y \sim_L w$ implies that $\pr_K(y) \sim_{L,K} \pr_K(w)$. By our induction hypothesis,  $\pr_K(y) \sim_{L,K} \pr_K(w)$ implies that $\rxi(\pr_K(y)) = \rxi(\pr_K(w))$. So if $b/a \les n-1$, then we are done by (\ref{eqn 01}). If not, then it remains to determine the membership of $t_n$ in $\rxi(y)$ and $\rxi(w)$.

So, suppose now that $b/a > n-1$. Then by Corollary 6.7 of \cite{BI}, we know that $\ell_t(y) = \ell_t(w)$. By Theorem 1 of \cite{GeckInd}, we have $\ell_t(\pr_K(y)) = \ell_t(\pr_K(w))$. It follows that $\ell_t(\rep_K(y)) = \ell_t(\rep_K(w))$. It remains to observe that $t_n \in \rxi(w)$ if and only if $\ell_t(\rep_K(w)) = 1$. \qed

\begin{corollary}
Let $(W,S,\wei)$ be a finite weighted Coxeter system. Then $\rxi$ is an invariant of its left cells.
\end{corollary}

Proof -- We use the classification of finite irreducible Coxeter groups to reduce to the case where $W$ is of type $B_n$, type $F_4$, or type $I_2(m)$ with $m$ even. The case of type $F_4$ was looked at in Lemma~\ref{F_4 prop}. The case of type $I_2(m)$ has been considered in Example 7.4 of \cite{BG}. For type $B_n$, if $b/a \in (0,1)$, then $\rxi = \rdes^{\raise  0.5pt \hbox{\scriptsize $\wei$}}$, and we are in the case of   Example 7.4 of \cite{BG} again. Finally, we appeal to Proposition~\ref{SPD invariant} to conclude the proof. \qed

\ms

From the definition of $\rxi$, we see that we can write \[W_n \, = \!\bigsqcup_{I\, \subseteq\, \Sxi} \{\qs w \in W_n \, : \, \rxi(w)=I \qs\},\]

\noi and so $\rxi$ partitions $W_n$ into up to $2^{2n-1}$ subsets. However, unlike with $\rdes$, some of these subsets will be empty.
 For instance, if $W_2$ is equipped with an asymptotic weight function, then $\rxi$ distinguishes all six of its left cells. On the other hand, there is no element $w \in W_2$ such that $\rxi(w)$ is equal to either $\{t,s_1\}$ or $\{t_2\}$. These observations motivate us to work out just how fine an invariant $\rxi$ is for the left cells of $W_n$. 

\begin{lemma}\label{asy spd}
Suppose $W_n$ is equipped with an asymptotic weight function. Then $\rxi$ partitions 
$W_n$ into exactly $2\cdot3^{n-1}$ non-empty subsets.
\end{lemma}

Proof -- We proceed by induction on $n$. For the case $n=2$, see 
Example 7.4 of \cite{BG}. 

\ms

Let $x \in X_K$ and $u \in W_K$. Then by Remark~\ref{incl}, $\rxi(xu)$ is equal to one of the following:
\[ \begin{array}{rlrl} (i) & \!\! \rxi(u),\qq  & (ii) & \!\!\rxi(u)\cup\{s_{n-1}\},\\  (iii) & \!\!\rxi(u)\cup\{t_n\},\qq & (iv) & \!\!\rxi(u)\cup\{s_{n-1},\,t_n\}.\end{array} \]


Let $I \subseteq \Sxi$, and consider the sets $C_I:=\{ u \in W_K \ : \ \rxi(u) = I\}$ and $D_I:=\{xu \ : \ x \in X_K,\ u \in C_I\}$. By our inductive hypothesis, $\rxi$ partitions $W_K$ into exactly $2 \cdot 3^{n-2}$ non-empty subsets. It suffices to show that if $C_I \neq \emptyset$, then $\rxi$ partitions $D_I$ into exactly three non-empty subsets. So, fix some $I \subseteq \Sxi$ such that $C_I \neq \emptyset$, and let $u \in C_I$.

\ms

First note that $e \in X_K$ is such that $s_{n-1},t_n \notin \rxi(e  u)$ for all $u \in W_K$, while $t_n \in X_K$ commutes with all $u \in W_K$ to give $s_{n-1},t_n \in \rxi(t_nu)$. We now distinguish two cases. 

\ms

Case A -- Suppose that $t_{n-1} \notin I$. By Lemma~\ref{signdes}, we have $u(n-1)>0$. Consider $x = s_1 \cdots s_{n-1} \in X_K$. 
We have $xu(n-1)>0$ and $xu(n) = 1$, and it follows that 
$\rxi(xu) = \rxi(u)\cup\{s_{n-1}\}$.

So now suppose there is some $x' \in X_K$ such that $s_{n-1} \notin \rxi(x'u)$. Using Lemma~\ref{signdes} and  Remark~\ref{incl}, we have that $x'u(n)>x'u(n-1)>0$, and so $t_n \notin \rxi(x'u)$. Thus if $t_{n-1} \notin I$, then $\rxi(x'u)$ is not equal to $\rxi(u) \cup \{t_n\}$.

\ms

Case B -- Suppose that $t_{n-1} \in I$, and let $x = t s_1 \cdots s_{n-1} \in X_K$. Via similar considerations to those in Case A, we see that $\rxi(xu) = \rxi(u)\cup\{t_n\}$ while $\rxi(x'u) \neq \rxi(u)\cup\{s_{n-1}\}$ for all $x' \in X_K$. \qed

\begin{corollary}
Suppose $b/a \in (k,k+1]\subseteq(1,n]$ for some $k \in \Z$. Then $\rxi$ partitions 
$W_n$ into exactly $2^{n-k}\cdot3^k$ non-empty subsets.
\end{corollary}



\section{$\Area_n$}\label{S5}

In this section, we collect a number of technical results for later use. 

\subsection{Shape}\label{S2.1}  To any $w \in W_n$ we associate the bipartition $\sh(w)$ as in \S\,\ref{B_n section}. For any bipartition $\lam \Vdash n$, we set $\Omega_\lam:=\{w \in W_n \,:\, \sh(w) = \lam\}$, and for $0 \les q \les n$ we set $\ze_q:=(1^{n-q} \,|\, 1^q) \Vdash n$. 
If we are in the asymptotic case, then $\Omega_\lam$ is equal to a two-sided cell, as in Theorem~\ref{ThmBI} (iii).

Set $J:=\{s_1, \ldots, s_{n-1}\}$, and denote by $w_J \in W_n$ the longest word of the parabolic subgroup $W_J$. Then $\sh(w_J) = \ze_0$ and $\sh(w_Jw_0) = \ze_n$. 
We note that if $w \in W_n$ is such that either $A_n(w)=A_n(w_J)$ or $B_n(w)=B_n(w_J)$, then $w=w_J$. By Theorem~\ref{ThmBI}, the element $w_J$ lies in an asymptotic  left, right and two-sided cell of cardinality one. 
Analogous statements hold for the element $w_Jw_0$.

\begin{definition}\label{area lr cells}
\[\begin{aligned}
\ \ \Area_n :=\, & \, {\scalebox{0.6}{$\!$}}\bigsqcup_{q{\scalebox{0.6}{$\,$}}={\scalebox{0.6}{$\,$}}0}^n \Omega_{\ze_q},  \\
\rArea_n:=\, & \,\Area_n \setminus\{\,w_J,\, w_Jw_0\,\}.
\end{aligned}\]
\end{definition}

We have $\Area_n = w_0 \Area_n$ and $\Area_n = \{ w^{-1} \,:\, w \in \Area_n\}$, as well as analogous statements for $\rArea_n$.

\ms

The cardinality of $\{ B_n(w) \,:\, w \in \Omega_{\ze_q}\}$ (and of $\{A_n(w) \,:\, w \in \Omega_{\ze_q}\}$) is~$\nchooseq$. If we are in the asymptotic case, this is equivalent to saying that~$\Omega_{\ze_q}$ contains 
$\nchooseq$ left (equivalently, right) cells, each of size 
$\nchooseq$. For later use, we may therefore state:
\begin{equation}
\label{area number of cells 2}  \ \ \big|\{B_n(w) \,:\, w \in \Area_n\}\big| \eq \sum_{q{\scalebox{0.6}{$\,$}}={\scalebox{0.6}{$\,$}}0}^{n} {n \choose q} \eq 2^n.
\end{equation}

\subsection{Signed permutations}\label{S2.3} The condition $w \in \Area_n$ is quite a strict one, and in turn places conditions on the row form of $w$. Indeed, let us consider 
$w = w(1), \ldots, w(n) \in \Area_n$ 
as a signed permutation. Then the subsequence $x_1, \ldots, x_q$ of negative integers 
must be such that $|x_1| > |x_2| > \cdots > |x_q|$. 
Similarly, the subsequence $y_1, \ldots, y_{n-q}$ of positive integers must be such that $y_1 > y_2 > \cdots > y_{n-q}$. Further, if $y, w \in \Area_n$ then $B_n(y) = B_n(w)$ if and only if $y(i)$ and $w(i)$ have the same sign for all $1 \les i \les n$.

\begin{example} We have $y,w \in \Area_7$, where: 
\[\! \begin{array}{rrl}
 y\qs =\left( \begin{array}{ccccccc} 1 & 2 & 3 & 4 & 5 & 6 & 7 \\ -7 & -5 & 6 &4 &3 & -2 & 1 \end{array} \right)\!, \!\! & 
\Yvcentermath2  \!\!A_7(y) =\  \qc\qc \young(1,3,4,6)\ \  {\raise -10pt \hbox{\young(2,5,7)}}\,, & 
   \Yvcentermath2 \!\! \qs B_7(y) =\ \qc\qc \young(3,4,5,7)\ \  {\raise -10pt \hbox{\young(1,2,6)}}\, ,  \\[2.3 em]
 w\qs =\left( \begin{array}{ccccccc} 1 & 2 & 3 & 4 & 5 & 6 & 7 \\ -4 & -2 & 7 &6 &5 & -1 & 2 \end{array} \right)\!, \! &
\Yvcentermath2  \!\!A_7(w) =\  \qc\qc \young(3,5,6,7)\ \  {\raise -10pt \hbox{\young(1,2,4)}}\,, & \!
   \Yvcentermath2 \!\! B_7(w) =\ \qc\qc \young(3,4,5,7)\ \  {\raise -10pt \hbox{\young(1,2,6)}}\,.
\end{array}\]
\end{example}

\subsection{Reduced expressions}\label{S2.4}
The forthcoming notation $\si_{n,q}$ and $p_{n,q}$ will be crucial in discussions regarding $\Area_n$. First set $a_0:=e$, $b_{n-1}:=b_n:=e$, and $p_{n,0}:=p_{n,n}:=e$, and then:
\begin{itemize}
\item for $1 \les q \les n$, set $a_q:= (t) (s_1t) \cdots (s_{q-1} \cdots s_1t)$, \\[-0.9em]
\item for $0 \les q \les n-2$, set $b_q:= (s_{q+1}) (s_{q+2} s_{q+1}) \cdots (s_{n-1} \cdots s_{q+1})$, \\[-0.9em]
\item for $0 \les q \les n$, set $\si_{n,q} := a_q \cdot b_q$, \\[-0.9em]
\item for $1 \les q \les n-1$, set:
\end{itemize}
\[ \begin{aligned} p_{n,q}:=\, &  \qs (s_{n-q} s_{n-q-1} \cdots s_1) (s_{n-q+1} s_{n-q} \cdots s_2) \cdots (s_{n-1} s_{n-2} \cdots s_q)  \\
  \eq &  \qs (s_{n-q} s_{n-q+1} \cdots s_{n-1}) (s_{n-q-1} s_{n-q} \cdots s_{n-2}) \cdots (s_1 \cdots s_q).\end{aligned}\]
 Note that $a_q$ and $b_q$ commute in $W_n$, with $\ell(\si_{n,q}) = \ell(a_q) + \ell(b_q)$. 

\begin{proposition}\label{asy cells in area}
Consider the asymptotic two-sided cell $\Omega_{\ze_q} \subseteq \Area_n$. It contains the asymptotic left cell
\[ \Ga_q:= \{\, \pi \si_{n,q} \,:\, \pi \les_e p_{n,q}\,\},\]
and the following is a complete list of asymptotic left cells contained in $\Omega_{\ze_q}$:
\[ \{\, \Ga_q \tau^{-1} \,:\, \tau \les_e p_{n,q}\,\}.\]
Thus,
\[ \Omega_{\ze_q} \eq \{\, \pi \si_{n,q} \tau^{-1} \,:\, \pi,\tau \les_e p_{n,q}\,\}.\]
Further, for any  $0 \les q \les n$ and $\pi,\tau \les_e p_{n,q}$, the expression  $\pi \si_{n,q} \tau^{-1}$ is reduced.
\end{proposition}

Together with \S\,\ref{S2.1}, this result implies that $p_{n,q}$ has 
$\nchooseq$ suffixes.


\subsection{Cell decomposition}\label{S2.5}
We now give a pair of basic lemmas with useful corollaries.

\begin{lemma}\label{TL}
Let $w \in \Area_n$. Then: \\[-0.7em]

\setlength{\tabcolsep}{4pt}
\hskip -5.5pt \begin{xtabular}{r p{319 pt}}

(i) &  If $wt > w$ and $wt \in \Area_n$, then $s_1 \in \rdes(w)$ and $\rdes(wt)=\rdes(w)\cup\{t\}\setminus\{s_1\}$.\\[0.1em]

(ii) &  For $1 \les i \les n-2$, if $ws_i > w$ and $ws_i \in \Area_n$, then $s_{i+1} \in \rdes(w)$. We have $\rdes(ws_1) = \rdes(w) \cup\{s_1\} \setminus\{t,s_2\}$, and for $2 \les i \les n-2$ we have  $\rdes(ws_i) = \rdes(w) \cup\{s_i\} \setminus\{s_{i+1}\}$. \\[0.1em]

(iii) &  If $ws_{n-1} > w$ and $ws_{n-1} \in \Area_n$, then $\rdes(ws_{n-1}) = \rdes(w) \cup\{s_{n-1}\}$. \\[0.1em]

(iv) &  Left-handed versions of the above statements also hold.
\end{xtabular}
\end{lemma}

Proof -- Following the discussion in \S\,\ref{S2.3}, this may be verified with judicious application of Lemma~\ref{signdes}. \qed

\begin{corollary}\label{TL cor}
Let $y,w \in \Area_n$ with $\rdes(y)=\rdes(w)$, and let $p \in W_n$. Suppose that
\[ y \tau^{-1},\ w\tau^{-1} \in \Area_n \ \ \ \ \ \ \ \forall\, \tau \les_e p.\]
Then we have $\rdes(y\tau^{-1}) \eq \rdes(w\tau^{-1})$ for all $\tau \les_e p$.
\end{corollary}


\begin{lemma}\label{TL2}
Let $w \in \Area_n$, and $1 \les i \les n-2$. Suppose that $s_iw \in \Area_n$ with $w < s_iw$. Then $w \sim_L s_iw$ with respect to any choice of parameters.
\end{lemma}

Proof -- Application of Lemma~\ref{TL} (iv) shows that $M^{s_{i+1}}_{w,s_iw} = 1$, and thus $w \sim_L s_iw$, with respect to any choice of parameters.\footnote{\ For the definition of $M$-polynomials, the reader may consult Chapter 6 of \cite{L}, where they are denoted by $\mu$.}\qed

\ms

Denote by $S_w$ the set of all $s \in S$ such that $s$ occurs in a reduced expression for $w$ (as in \S\,9.2 of \cite{L}).

\begin{corollary}\label{subcells}
Let $\Ga \subseteq \Area_n$ be an asymptotic left cell. Then the cell can be decomposed as
\[\Ga \,=\, \ga_1 \sqcup \ga_2\]
where $\ga_1$ and $\ga_2$ are each contained within a left cell for all choices of parameters.
Denote by $\si$ the element in $\Ga$ of minimal length, and set $q=\ell_t(\si)$.
Then we have
\[ \begin{array}{l}
\ga_1 \eq \{\, \pi \si \,:\, \pi \les_e p_{n,q} \,\text{ and }\, s_{n-1} \notin S_\pi\,\} \eq \{\,\pi \si \,:\, \pi \les_e p_{n-1,q}\,\}, \\
\ga_2 \eq \{\, \pi \si \,:\, \pi \les_e p_{n,q} \,\text{ and }\, s_{n-1} \in S_\pi\,\} \eq \{\,\pi \chi_q \si \,:\, \pi \les_e p_{n-1,q-1}\,\},
\end{array}\]
where $\chi_q:= s_{n-1} \cdots s_q$.
\end{corollary}

Proof -- We first check that the two descriptions of $\ga_1$ and $\ga_2$ are consistent. By Proposition~\ref{asy cells in area}, we may describe the cell $\Ga$ as
\[ \Ga \eq \{\, \pi \si \,:\, \pi \les_e p_{n,q}\,\}.\]
Recall the following reduced expressions for $p_{n,q}$:
\[\begin{aligned}
p_{n,q} & \eq (s_{n-q} \cdots s_{n-1})  (s_{n-q-1} \cdots s_{n-2}) \cdots (s_1 \cdots s_q) \\
& \eq (s_{n-q} \cdots s_1)  (s_{n-q+1} \cdots s_2)\cdots (s_{n-1} \cdots s_q).
\end{aligned}\]
In the first expression, we have moved the lone $s_{n-1}$ term as far to the left as possible; in the second, as far to the right as possible. We see that $s_{n-1} \notin S_\pi$ if and only if $\pi \les_e p_{n-1,q}$. Similarly, $s_{n-1} \in S_\pi$ if and only if there exists a reduced expression for $\pi$ ending in the block sequence $\chi_q = s_{n-1} \cdots s_q$. Now note that
\[ p_{n,q} = p_{n-1,q-1} \chi_q \ \text{ with } \ \ell(p_{n,q}) \eq \ell(p_{n-1,q-1}) + \ell(\chi_q).\]

Following this, we may apply Lemma~\ref{TL2} to see that for $i \in \{1,2\}$, if $y,w \in \ga_i$, then $y \sim_L w$ for any choice of parameters. 
In this result, we note that $\ga_2 = \emptyset$ if $|\Ga|=1$. \qed



\subsection{The partition of $\Area_n$ with respect to $\rxi$}\label{S2.6}

\begin{lemma}\label{s and t relation}
Let $w \in \Area_n$, and let $1 \les i \les n-1$. Then:
\[ \ell(ws_i) < \ell(w) \iff \ell(wt_i) > \ell(w).\]
\end{lemma}

Proof -- The statement follows from Lemma~\ref{signdes} and \S\,\ref{S2.3}. \qed

\begin{corollary}\label{corollary 5.9}
Consider $(W_n,S,\wei)$ and let $y,w \in \Area_n$. 

\setlength{\tabcolsep}{4pt}
\hskip -5.5pt \begin{tabular}{r p{319pt}}
(i) &  If $b/a > n-2$, then
\[ \sm\sm\sm\sm\sm\sm\sm\sm\!   \rxi(y) \cap \{t_1, \ldots, t_n\} = \rxi(w) \cap \{t_1, \ldots, t_n\} \iff \rxi(y)=\rxi(w).\]  \\[-1em]
(ii) & If $b/a \in [1,n-1]$, then
\[ \sm\sm\sm\sm\sm\sm\sm\sm\sm\sm\!  \rxi(y)=\rxi(w) \iff \rdes(y)=\rdes(w).\]
\end{tabular}
\end{corollary}


\begin{lemma}\label{lemma 5.10}
Consider $(W_n,S,\wei)$.

\setlength{\tabcolsep}{4pt}
\hskip -5.5pt \begin{tabular}{r p{319pt}}
(i) & If $b/a > n-1$, then $\Area_n$ is partitioned by $\rxi$ into $2^n$ non-empty subsets; these subsets are exactly the asymptotic left cells. \\

(ii) & Suppose $b/a > n-1$, and let $T \subseteq \{t_1, \ldots, t_n\}$ be arbitrary. Then there exists a unique left cell $\Ga\subseteq \Area_n$ such that for any $w \in \Ga$, we have $\rxi(w) \cap \{t_1, \ldots, t_n\}=T$. \\

(iii) & If $b/a \in [1,n-1]$, then $\Area_n$ is partitioned by $\rxi$ (and $\rdes$) into $2^{n-1}$ non-empty subsets, each of which is the union of exactly two asymptotic left cells. 
\end{tabular}

\end{lemma}

Proof -- We first look at part (i). By Corollary~\ref{corollary 5.9} (i), we know that 
$\rxi(y) =\rxi(w)$ if and only if $\rxi(y) \cap \{t_1, \ldots, t_n\} = \rxi(w) \cap \{t_1, \ldots, t_n\}$. 
By Lemma~\ref{signdes}, this is in turn equivalent to the statement that $y(i)$ and $w(i)$ have the same sign for all $1 \les i \les n$. As was noted in \S\,\ref{S2.3}, this is if and only if $B_n(y)=B_n(w)$, which is true if and only if $y$ and $w$ are in the same asymptotic left cell. Recalling (\ref{area number of cells 2}) from \S\,\ref{S2.1} finishes the proof.

Part (ii) follows from part (i), Corollary~\ref{corollary 5.9} (i), and the fact that the power set of $\{t_1, \ldots, t_n\}$ contains $2^n$ elements.

Turning to part (iii), we first prove the statement for $b/a \in (n-2,n-1]$. So, let $T' \subseteq \{t_1, \ldots, t_{n-1}\}$ be an arbitrary subset. Then
\[  \begin{aligned} & \{\,w \in \Area_n \,:\, \rxi(w) \cap\{t_1, \ldots, t_{n-1}\} = T'\,\} \text{ (with $b/a \in (n-2,n-1]$) }  \\
= \ & \{\,w \in \Area_n \,:\, \rxi(w) \cap\{t_1, \ldots, t_n\} = T'\,\} \\
& \sqcup\ \{\,w \in \Area_n \,:\, \rxi(w) \cap\{t_1, \ldots, t_n\} = T'\sqcup\{t_n\}\,\} \text{ (with $b/a > n-1$).}
\end{aligned}\]
The parts of this disjoint union are non-empty by part (ii). So when $b/a \in (n-2,n-1]$, $\rxi$ partitions $\Area_n$ into half as many non-empty subsets compared to when $b/a > n-1$. Now apply part (i). Corollary~\ref{corollary 5.9}~(ii) then extends the scope of the result to all $b/a \in [1,n-1]$.  \qed

\ms

Later, it will benefit us to have a description of how $\rxi$ partitions $\Area_n$. Lemma~\ref{lemma 5.10} indicates that we can split this into two cases. If $b/a > n-1$, then descriptions of this partition are given in both Theorem~\ref{ThmBI}~(i) and Proposition~\ref{asy cells in area}. Otherwise, it suffices to look at how $\rdes$ partitions $\Area_n$. 

 So let $I \subseteq S$ be such that $\{\qs w \in \Area_n \,:\, \rdes(w) = I \qs\}$ is non-empty. 
 This set comprises exactly two asymptotic left cells. In order to describe all such sets in this form, we need to show how the left cells `pair up' under~$\rdes$. 
This has three steps; first finding representative pairings, then taking advantage of Corollary~\ref{TL cor} to `translate' the resulting sets (pairs of cells) via right multiplication,  and finally showing that this method encompases all elements in $\Area_n$. Some notation:

\begin{itemize}
\item denote by $\Ga(w)$ the asymptotic left cell containing $w \in W_n$, \\[-1em]
\item denote by $\Ga'(w)$  the asymptotic right cell containing $w \in W_n$, \\[-1em]
\item set $\Up(w):=\{ z \in \Area_n \,:\, \rdes(z)=\rdes(w)\}$ for $w \in \Area_n$,  \\[-1.1em]
\item and recall that $\chi_q := s_{n-1} \cdots s_q$.
\end{itemize}

Let us begin by using Lemma~\ref{signdes} to see that $\rdes(\si_{n,q}) = \{t,s_{q+1},\ldots,s_{n-1}\}$. As $\chi_{q+1} \les_e p_{q+1,n}$,  Proposition~\ref{asy cells in area} indicates that $\si_{n,q+1}\chi_{q+1}^{-1} \in \Area_n$, and by applying Lemma~\ref{signdes} again we see that $\rdes(\si_{n,q}) = \rdes(\si_{n,q+1}\chi_{q+1}^{-1})$. 

\ms

Thus for $0 \les q \les n-1$, we have:
\[ \Up(\si_{n,q}) \eq \Ga(\si_{n,q}) \,\sqcup\, \Ga(\si_{n,q+1}\chi^{-1}_{q+1}).\]
By Proposition~\ref{asy cells in area}, we know that $\Ga'(\si_{n,q}) = \{\, \si_{n,q} \pi^{-1} \,:\, \pi \les_e p_{n,q}\, \}$.
As $p_{n-1,q} \les_e p_{n,q}$ and $\Ga'(\si_{n,q}) \subseteq \Area_n$, we observe that:
\[ \{\, \si_{n,q} \tau^{-1} \,:\, \tau \les_e p_{n-1,q}\,\} \subseteq \Area_n.\]
Similarly, we have $\Ga'(\si_{n,q+1}) \eq \{\, \si_{n,q+1} \pi^{-1} \,:\, \pi \les_e p_{n,q+1} \,\}$.
Noting that $p_{n,q+1} = p_{n-1,q} \chi_{q+1}$ with $\ell(p_{n,q+1}) = \ell(p_{n-1,q}) + \ell(\chi_{q+1})$, we also observe that:
\[ \{\,\si_{n,q+1} \chi^{-1}_{q+1} \tau^{-1} \,:\, \tau \les_e p_{n-1,q} \,\} \subseteq \Area_n.\]

These two observations allow us to apply Corollary~\ref{TL cor} (setting $y=\si_{n,q}$, $w= \si_{n,q+1}\chi^{-1}_{q+1}$ and $p=p_{n-1,q}$) to see that for all $0 \les q \les n-1$ and all $\tau \les_e p_{n-1,q}$ we have:
\[\Up(\si_{n,q} \tau^{-1}) \eq \Ga(\si_{n,q} \tau^{-1}) \,\sqcup\, \Ga(\si_{n,q+1} \chi^{-1}_{q+1} \tau^{-1}).\]

Sets procured in this way are non-empty, mutually distinct and contained in $\Area_n$. Recalling from \S\,5.4 that the number of suffixes of $p_{n,q}$ is 
$\nchooseq$, we see that the number of sets that we have obtained is $2^{n-1}$, 
which by Lemma~\ref{lemma 5.10} (ii) is the number of non-empty sets that $\rdes$ partitions $\Area_n$ into. 
We therefore obtain Proposition~\ref{intermediate area cells} parts (i) and (ii) below. 

There is a more straightforward description of these sets.
Using the braid relations for $W_n$, we can verify that for $0 \les q \les n-1$, we have
\[ s_{n-q-1} \cdots s_1 \cdot t \cdot p_{n,q}  \eq p_{n,q+1} \cdot t \cdot s_1 \cdots s_q.\]
Denote this element by $\Pi_{n,q}$, and note that both of these expressions are reduced.

\begin{proposition}\label{intermediate area cells}
Let $w \in \Area_n$.  \\[-0.5em]

\setlength{\tabcolsep}{4pt}
\hskip -5.5pt \begin{xtabular}{rl}
(i) & There exists $0 \les q \les n-1$ and $\tau \les_e p_{n-1,q}$ such that: \\[0.8em] \noi & 
$ \displaystyle \Up(w) \eq \Up(\si_{n,q} \tau^{-1}) \eq \Ga(\si_{n,q} \tau^{-1}) \,\sqcup\, \Ga(\si_{n,q+1} \chi^{-1}_{q+1} \tau^{-1})$, \\ [0.5em]

(ii) & $ \displaystyle  \Area_n \eq \bigsqcup_{q{\scalebox{0.6}{$\,$}}={\scalebox{0.6}{$\,$}}0}^{n-1}\ \bigsqcup_{\tau\, \les_e\, p_{n-1,q}} \sm \Up(\si_{n,q} \tau^{-1})$. \\[2em]

& For $0 \les q \les n-1$ and $\tau \les_e p_{n-1,q}$, we have: \\[0.5em]

(iii) & $ \displaystyle \Up(\si_{n,q} \tau^{-1}) \eq \{\, \pi \si_{n,q} \tau^{-1} \,:\, \pi \les_e \Pi_{n,q} \,\}$, \\ [0.5em]

(iv) & $ \displaystyle |\Up(\si_{n,q} \tau^{-1})| \eq {n \choose q} + {n \choose q+1}$. \\[0.5em]

\end{xtabular} 

\end{proposition}

Proof -- Verifying that the descriptions of $\Up(\si_{n,q} \tau^{-1})$ in parts (i) and (iii) are equivalent is done in an entirely similar way to verifying the asymptotic left cell decomposition in the proof of Corollary~\ref{subcells}. For part (iv), consider the decomposition of $\Up(\si_{n,q} \tau^{-1})$  in part (i) into asymptotic cells. 
\qed


\section{An extension of the generalised $\tau$-invariant in type $B_n$}
A classical pair of results in the theory of Kazhdan--Lusztig cells with equal parameters are as follows: if two elements of a Coxeter group lie in the same left cell then they have the same generalised $\tau$-invariant, and two such elements remain equivalent under $\sim_L$ after a $*$-operation has been applied to them; see \S\,3 of Vogan \cite{Vogan} or \S\,4 of  \cite{KL}. The theory of $*$-operations and the generalised $\tau$-invariant has been substantially generalised in \cite{BG}, providing compatibility with unequal parameters as well as the framework for maps richer than the traditional $*$-operations. In this section, we recall some definitions and results from \cite{BG}, and introduce new ones. 

\subsection{Vogan classes}\label{S3.1}

Let $(W,S,\wei)$ be a weighted Coxeter system. 

\begin{definition} (Bonnaf\'{e}--Geck,  Definition 6.1 in \cite{BG}) A pair $(I, \de)$ consisting of a non-empty subset $I\subseteq S$ and a left cellular map $\de:W_I\ra W_I$ is called {\em KL-admissible}. We recall that this means that the following conditions are satisfied for every left cell $\Ga \subseteq W_I$ (with respect to $\wei|_{W_I}$):

\ms

(A1) \qs $\de(\Ga)$ also is a left cell.

(A2) \qs The map $\de$ induces an $\IH_I$-module isomorphism\footnote{\ We denote by $[\Ga]$ the $\IH$-module afforded by a left cell $\Ga \subseteq W$; see \S\,6 of \cite{LWeyl} for details.} 
 $[\Ga] \cong [\de(\Ga)]$.

\ms

\noi We say that $(I,\de)$ is {\em strongly KL-admissible} if, in addition to (A1) and (A2), the following condition is satisfied:

\ms

(A3) \qs We have $u \sim_{R,I} \de(u)$ for all $u \in W_I$.
\end{definition}

\noi If $I \subseteq S$ and if $\de:W_I \ra W_I$ is a map, we obtain a map $\de^L:W \ra W$ by
\[ \de^L(xu) := x \de(u) \ \ \ \ \text{for all $x \in X_I$ and $u \in W_I$.}\]
The map $\de^L$ is called the {\em left extension of $\de$ to $W$}. However, by abuse of notation, we will often use $\de$ to refer to $\de^L$ where the meaning is clear.

\begin{theorem}\label{Thm BG 6.2} (Bonnaf\'e--Geck, Theorem 6.2 in \cite{BG})
Let $(I,\de)$ be a (strongly) KL-admissible pair. Then $(S,\de^L)$ is (strongly) KL-admissible.

\end{theorem}

This theorem brings us to consider strongly KL-admissible pairs that give us as much information about the cells of $W$ as possible; to this end, we introduce an additional condition.

\ms

(A4) \qs If $u,v \in W_I$ are such that $u \sim_{R,I} v$, then there exists
some $k \in \Z_{\ges0}$ such that $u = \de^k(v)$.

\begin{definition}
We say that a pair $(I,\de)$ is {\em maximally KL-admissible} if conditions (A1), (A2), (A3) and (A4) hold.
\end{definition}

\begin{example}\label{ex S_n}
 Consider $W_n$ equipped with any weight function, and let $J:=\{s_1,\ldots, s_{n-1}\}$, so that $W_J \cong \mathfrak{S}_n$. The cells of $\mathfrak{S}_n$ are described by the Robinson--Schensted correspondence; see \S\,5 of \cite{KL} or Ariki \cite{Ariki}.

Let $\lambda\vdash n$ be a partition, and $\Omega_\lambda$ be the two-sided cell of $W_J$ such that $\operatorname{sh}(w) = \lambda$ for any $w \in \Omega_\lambda$. Fix an arbitrary total order on the left cells contained in $\Omega_\lambda$; $\Ga_1 < \Ga_2 < \cdots < \Ga_k$.

We now define a map $\ep_\lam:\Omega_\lam \ra \Omega_\lam$. Let $1 \les i  \les k-1$. For any $w \in \Ga_i$, set $\ep_\lam(w)$ to be equal to the unique element $w'  \in \Ga_{i+1}$ such that $w \sim_{R,J} w'$. For $w \in \Ga_k$, set $\ep_\lam(w)$ to be equal to the unique element $w' \in \Ga_1$ such that $w \sim_{R,J} w'$, as in Figure~\ref{maximally KL-admissible}.

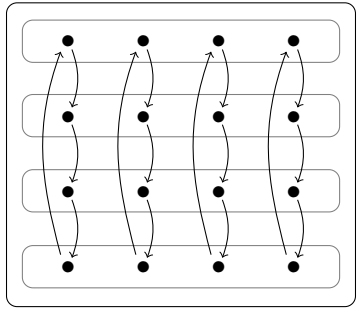
\begin{figure}
\[\begin{tikzpicture}
\tikzstyle{lrcell} = [rectangle, rounded corners, minimum width=132pt, minimum height=115pt, text centered, draw=black, fill=white!25]
\tikzstyle{lcell} = [rectangle, rounded corners, minimum width=120pt, minimum height=16pt,text centered, draw=gray ]

\node at (0,0) [lrcell] {};

\node at (0,-1.49) [lcell] {};
\node at (0,-.48) [lcell] {};
\node at (0,1.51) [lcell] {};
\node at (0, .52) [lcell] {};

\foreach \x in {0,1,2,3}
\foreach \y in {0,1,2,3}
{   \node at (\x - 1.5,\y -1.5) {$\bullet$};   }

\foreach \x in {0,1,2,3}
{   \arghk{\x-1.4 }{1.05}  }
\end{tikzpicture}\]

\caption{\small A two-sided cell $\Omega_\lam \subseteq \mathfrak{S}_n$ 
with left cells given by rows, and right cells by columns. The map 
$\ep_\lam$ is indicated by the arrows.}\label{maximally KL-admissible}
\end{figure}

This map can be then extended to the rest of the group by setting:
\[ \ep'_\lam(w) = \left\{ \begin{array}{rr} \ep_\lam(w) & \text{ if  } w \in \Omega_\lam, \\ w\phantom{)} & \text{ otherwise.} \end{array}\right.\]

Finally, we can define a map $\ep:W_J \ra W_J$ by taking the composition of all maps $\ep'_\lam:W_J \ra W_J$ over the indexing set $\{\lambda \ : \ \lambda \vdash n\}$.

Example 2.6 of Geck \cite{GeckRel} tells us that $\ep$ satisfies condition (A2). From the construction, we can see that $(J,\ep)$ is maximally KL-admissible, and well-defined up to the choice of total ordering on left cells in each two-sided cell.
\end{example}

\begin{example}\label{ex W_n}
Consider $(W_n,S,\wei)$ and set $K:=\{t,s_1, \ldots, s_{n-2}\}$ so that $W_K = W_{n-1}$.
 Suppose that the restriction of $\wei$ to $W_K$ is an asymptotic weight function for $W_K$. Then the left, right and two-sided cells of $W_K$ are determined by the generalised Robinson--Schensted correspondence, as in Theorem~\ref{ThmBI}. Now we can proceed analogously to Example~\ref{ex S_n} to obtain a map $\psi:W_K \ra W_K$. Theorem 6.3 of \cite{GeckRel} tells us that $\psi$ satisfies condition (A2), as does Example 6.8 of \cite{BG}. The resulting pair $(K, \psi)$ is maximally KL-admissible.

Recall from Theorem~\ref{ThmBI} that $b/a > n-1$ is the optimal lower bound for a weight function on $W_n$ to result in asymptotic cells. Thus, the maximally KL-admissible pair $(K,\psi)$ is well-defined (up to a choice of total ordering on left cells in two-sided cells) with respect to any $\wei:W_n \ra \Z_{\ges0}$ such that $b/a > n-2$.

Conversely, if $b/a \les n-2$, then $\wei|_{W_K}$ is not an asymptotic weight function for $W_K$. Without a suitable analogue of Theorem 6.3 of \cite{GeckRel}, we cannot readily produce a non-trivial KL-admissible pair $(\psi',K)$ with respect to $\wei$.
\end{example}

From now until the end of this section we fix an arbitrary weighted Coxeter system $(W,S,\wei)$, as well as a map $\rho:W \ra E$ (where $E$ is a fixed set) such that the fibres of $\rho$ are (possibly empty)\footnote{\ In \S\,7 of \cite{BG}, the condition that $\rho$ is surjective is not necessary.} unions of left cells.

\begin{definition}\label{def vogan} (Bonnaf\'e--Geck, \S\,7 of \cite{BG}) Let $\De$ be a collection of KL-admissible pairs with respect to $(W,S,\wei)$.
We define by induction on $n$ a family of equivalence relations $\app_n^{\De,\qs\rho}$ on $W$ as follows. Let $y,w \in W$.
\begin{itemize}
\item For $n=0$, we write $y \app_0^{\De,\qs\rho} w$ if $\rho(y)=\rho(w)$. \\[-0.8 em]
\item For $n \ges 1$, we write $y \app_n^{\De,\qs\rho} w$ if $y \app_{n-1}^{\De,\qs\rho} w$ and $\de^L(y) \app_{n-1}^{\De,\qs\rho} \de^L(w)$ for all $(I,\de) \in \De$.
\end{itemize}
We write $y \app^{\De,\qs\rho}w$ if $y \app_n^{\De,\qs\rho} w$ for all $n \ges 0$. The equivalence classes under this relation are called the {\em left Vogan} ($\De$, $\rho$){\em -classes}.
\end{definition}

\begin{theorem}\label{Thm BG 7.2} (Bonnaf\'e--Geck, Theorem 7.2 of \cite{BG})
Let $y,w \in W$. Then \[y \sim_L w \ \Ra \ y \approx^{\De,\,\rho}w.\]
\end{theorem}

\subsection{Orbits}\label{S3.2}
In the case of equal parameters, the $*$-orbit of an element $w$ is the set consisting of $w$ and all elements that can be reached from $w$ by a sequence of $*$-operations. As shown in \cite{KL}, the $*$-orbit of $w$ is contained in the right cell containing $w$.
Theorem~\ref{Thm BG 6.2} offers an analogue of this result in the more general framework of strongly KL-admissible pairs.

\ms

Let $\De$ be a collection of strongly KL-admissible pairs, and $y,w \in W$. If $y$ can be reached from $w$ via a sequence of applications of maps $\de^L$ such that $(I,\de) \in \De$, then we write $y \stackclose\De\llra w$, and we denote by $\orb^R_\De(w)$ the set of all elements related to $w$ in this way. We refer to this set as the {\em $\De$-orbit of $w$}, and it is a subset of the right cell containing $w$. By Theorem~\ref{Thm BG 6.2}  and Theorem~\ref{Thm BG 7.2} respectively, we have:
\[y^{-1} \stackrel{\De}\llra w^{-1} \ \Ra \ y \sim_L w \ \Ra \ y \approx^{\De,\,\rho}w.\]
Hence  left Vogan $(\De,\rho)$-classes are unions of left cells, which are in turn unions of inverses of $\De$-orbits.

\ms

Since inverses of  $\De$-orbits play an important role later, we also define:
\[\orb^L_\De(w)\,:=\,\{\,y \in W \,:\, y^{-1} \stackrel\De\llra w^{-1}\,\} \eq\{\,y^{-1} \in W \,:\, y \in \orb^R_\De(w^{-1})\,\}.\]
The superscripts are chosen so that $\orb^R_\De(w)$ is contained in the {\em right} cell containing $w$, and $\orb^L_\De(w)$ is contained in the {\em left} cell containing $w$.

\ms

Suppose that every $(I,\de) \in \De$ is such that the map $\de$ has finite order. It follows that $\de$ is a bijection and has an inverse. The left extension $\de^L$ shares these properties, so
\[ \mathcal{V}_\De \, := \, \langle\, \de^L \, :\, (I,\de) \in \De \,\rangle \]
is a group of permutations of the elements of $W$, and the orbit of $w \in W$ with respect to $\mathcal{V}_\De$ is precisely the $\De$-orbit of $w$.

\subsection{KL-admissible pairs}
The following lemma lists some properties of KL-admissible pairs. 
\begin{lemma}\label{KL pair lemma}
All sets of KL-admissible pairs mentioned will be with respect to $(W,S,\wei)$. 

\smallskip

\setlength{\tabcolsep}{4pt}
\hskip -5.5pt \begin{xtabular}{r p{319pt}}

(i)  &  If $(I,\de)$ is a strongly KL-admissible pair, then $\De := \{ (I,\de) \}$ may be replaced by a collection $\De'$ of strongly KL-admissible pairs, where each $(I',\de') \in \De'$ is such that $W_{I'}$ is an irreducible Coxeter group, and for all $y,w \in W$ we have:
\begin{itemize}
\item $y \stackrel\De\llra w \ \Ra \ y \stackrel{\,\qs \De'} \llra w$ \vskip 4pt
\item $y \app^{\De'\qc\qc,\,\rho} w \ \Ra \ y \app^{\De,\,\rho} w$.
\end{itemize} \\[-0.6em]

(ii) & Suppose that $(I,\de)$ is strongly KL-admissible and $(I',\de')$ is maximally KL-admissible, with $I \subseteq I'$. Then for all $y,w \in W$, we have: 
\begin{itemize}
\item $y \stackrel{\{(I,\qs\de)\}}\llra w \ \Ra \ y \stackrel{\{(I'\!,\qs\de')\}} \llra w$  \vskip 4pt
\item $y \app^{{\{(I'\!,\qs\de')\}},\,\rho} w \ \Ra \ y \app^{{\{(I,\qs\de)\}},\,\rho} w$.
\end{itemize} \\[-0.6em]

(iii) & Suppose that $W$ has a complete list $W_{I_1}, \ldots, W_{I_d}$ of distinct irreducible parabolic subgroups of rank $|S| - 1$, and for each $W_{I_i}$ there is a corresponding maximally KL-admissible pair $(I_i,\de_i)$. Set $\Xi:=\{ (I_i,\de_i) \,:\, 1 \les i \les d \qs\}$, and let $\De$ be any other pair of strongly KL-admissible pairs (such that\footnotemark\qs \qs for any $(I,\de) \in \De$ we have $I \neq S$). Then for all $w \in W$:
\begin{itemize}
\item $y \stackrel\De\llra w \ \Ra \ y \stackrel\Xi\llra w$ \vskip 4pt
\item $y \app^{\Xi,\,\rho} w \ \Ra \ y \app^{\De,\,\rho} w$.
\end{itemize}
\end{xtabular}
\end{lemma}
\footnotetext{\ Vogan classes are a tool for inductively obtaining cellular information from parabolic subgroups -- it makes sense, therefore, to restrict ourselves with the condition $I \neq S$.}

\vskip -15pt Proof -- We begin by noting that in all three parts, the first claim implies the second by Definition~\ref{def vogan}. We prove only part (iii).

Assume towards a contradiction that we have two elements $y,w \in W$ in the same $\De$-orbit but not in the same $\Xi$-orbit. So there is some $(I,\de) \in \De$ such that we have $y \in \orb^R_{\{(I,\de)\}}(w)$ but not $y \in \orb^R_{\Xi}(w)$. 
We have $y = (\de^L)^k(w)$ for some $k \in \Z_{\ges0}$; by the definition of the left extension of a map we must have $\rep_I(y) = \rep_I(w)$. Using part (i), it will suffice if we assume that $W_I$ is irreducible, and thus lies in some parabolic subgroup of rank $|S|-1$. So consider some $(I',\de') \in \Xi$ with $I \subseteq I'$, and apply part (ii) to obtain the contradiction.   \qed

\vfill

\subsection{Notation}\label{3.notation}\label{S3.3}
For the rest of this paper, we shall only work with weight functions $\wei:W_n \ra \Z\, $ 
 such that $b/a > n-2$. We also establish the following notation:
\begin{itemize}
\item let $(J,\ep)$ be as in Example~\ref{ex S_n}, \\[-1.1em]
\item let $(K,\psi)$ be as in Example~\ref{ex W_n},\\[-1.1em]
\item and set $\Xi:= \{ (J,\ep), (K,\psi) \}$.
\end{itemize}

The pair $(J,\ep)$ does not depend on our choice of weight function (the Iwahori--Hecke algebra of symmetric group does not admit unequal parameters), while we have mentioned in Example~\ref{ex W_n} that the pair $(K,\psi)$ is well-defined (up to a choice of ordering on the left cells) for any $\wei$ such that $b/a > n-2$. 

\ms

A corollary to Lemma~\ref{KL pair lemma} (ii) is that the $\Xi$-orbits  (and thus, left Vogan $(\Xi,\rho)$-classes) are independent of the choice of total ordering placed on the left cells in Examples~\ref{ex S_n}~and~\ref{ex W_n}. 

\ms

Noting that $\Xi$ satisfies the conditions in Lemma~\ref{KL pair lemma} (iii) motivates this particular choice of set of KL-admissible pairs. 

\ms

As the maps $\ep$ and $\psi$ have finite order, $\mathcal{V}_\Xi$ is a group of permutations of the elements of $W_n$.

\ms

Finally, we note that for the rest of this paper we will take $\rho = \rxi$.

\vfill


\section{Determining $\Xi$-orbits}\label{S7}

This section will be used to describe the decomposition of $W_n$ into $\Xi$-orbits. The key to this is relating the group $\mathcal{V}_\Xi$ to (generalised) Knuth relations, and thus to (bi)tableaux and cells. For $w \in \mathfrak{S}_n$, denote by $(P(w),Q(w))$ the application of the Robinson--Schensted algorithm.

\subsection{Generalised Knuth relations}
Let $w = w(1), \ldots, w(n) \in W_n$ be a signed permutation. We define the following relations, being modified versions of the ones found in 5.1.4 of Knuth \cite{Knuth} and \S\,3 of \cite{BI}.

\ms

For $1 \les i \les k-2$: if $w(i+1) < w(i) < w(i+2)$  or $w(i+2) < w(i) < w(i+1)$, then set $w':=ws_{i+1}$
and say that $w$ and $w'$ differ by a relation of type $\operatorname{I}_k$.

\ms

For $1 \les i \les k-2$: if $w(i+1) < w(i+2) < w(i)$ or $w(i) < w(i+2) < w(i+1)$, then set $w':=ws_{i}$ 
and say that $w$ and $w'$ differ by a relation of type $\operatorname{II}_k$.

\ms

For $1 \les i \les k-1$: if the sign of $w(i)$ and $w(i+1)$ differ, then set $w':=ws_{i}$ 
and say that $w$ and $w'$ differ by a relation of type $\operatorname{III}_k$.


\newpage

If two (signed) permutations $y,w$ are linked by a series of Knuth relations of, say, types I$_k$ and II$_k$, then we denote this by $y \sim_{\operatorname{I}_k,\,\operatorname{II}_k }w$. We may now state two important results in this context.

\begin{theorem}\label{K Knuth} (Knuth, Theorem 6 of \cite{Knuth})
Let $u,v \in W(A_{n-1})$. Then $u \sim_{\operatorname{I}_n,\,\operatorname{II}_n }v$ if and only if $P(u)=P(v)$.
\end{theorem}

\begin{proposition}\label{BI Knuth} (Bonnaf\'e--Iancu, Proposition 3.8 of \cite{BI})
Let $y,w \in W_n$. Then  $y \sim_{\operatorname{I}_n,\,\operatorname{II}_n ,\,\operatorname{III}_n}w$ if and only if $A_n(y)=A_n(w)$.
\end{proposition}

We will determine to what degree we can recover the Knuth relations using only the maps $\ep^L$ and $\psi^L$.

\begin{lemma}\label{type A proj lemma}
Let $w \in W_n$, and denote $u:=\pr_J(w)$. Then $w(i) > w(j)$ if and only if $u(i) > u(j)$, for $1 \les i, j \les n$.
\end{lemma}

Proof  -- We identify $x:=\rep_J(w)$. Let $p_1,\, p_2, \ldots,\, p_q$ be the increasing sequence comprising the $q=\ell_t(w)$ negative entries in the row form of~$w$, and similarly $m_1,\, m_2, \ldots,\, m_{n-q}$ the increasing sequence comprising the positive entries. Set $r_1:=t$ and recursively define $r_{i+1}:=s_ir_i=s_i \cdots s_1t$. Then from \S\,4.1 of \cite{BI}, we have
\[ x \eq r_{|p_q|} \!\cdots r_{|p_1|} \,=\, \left( \begin{array}{ccccccccc} 1 & 2 &  \cdots & q & q+1 & \cdots & n-1 & n \\[0.1em]
p_1 & p_2 & \cdots & p_q & m_{1} & \cdots & m_{n-q-1} & m_{n-q}
\end{array}\right).\]

\noi We can read off from this that the element $x$ has the property of preserving the ordering of the row form of $u$ under left multiplication; that is to say $u(i) > u(j)$ if and only if $xu(i) > xu(j)$. 
\qed

\begin{lemma}\label{type B proj lemma}
Let $w \in W_n$ and denote $u:=\pr_K(w)$.
Then

\smallskip

\setlength{\tabcolsep}{4pt}
\hskip -5.5pt \begin{tabular}{rl}
(i) & $w(i) > w(j)$ if and only if $u(i) > u(j)$,  for $1 \les i,j \les n-1$, \\[0.15em]

(ii) & $w(j) > 0$ if and only if $u(j) > 0$, for $1 \les j \les n-1$.
\end{tabular}

\end{lemma}

Proof -- Let $k = w(n)$. There exists a unique $x \in X_K$ such that $x(n) = k$, and so $x = \rep_K(w)$. Denote by $\sgn(z)$ the sign of $z \in \Z$. Then we have
\[ x \,=\, \left( \begin{array}{cccccccc} 1 & 2 & \cdots & |k|-1 & |k| &  \cdots & n-1 & n \\
1 & 2 & \cdots & |k|-1 & |k| + 1 & \cdots & n & k
\end{array}\right),\]
and in turn,
\[ w(i) \,=\, \left\{ \begin{array}{ll}
u(i) &  \text{ if } 1 \les |u(i)| \les |k|-1,  \\[0.1em]
u(i) + \sgn(u(i)) & \text{ if } |k| \les |u(i)| \les n-1, \\[0.1em]
k & \text{ if } i= n.
\end{array}\right.\]

\noi  Parts (i) and (ii) follow from this formula and Lemma~\ref{signdes}.   \qed

\begin{lemma}\label{lemma rep relations}
Let $y,w \in W_n$.

\smallskip

\setlength{\tabcolsep}{4pt}
\hskip -5.5pt \begin{tabular}{rl}

(i) & If $y \sim_{\operatorname{I}_n,\,\operatorname{II}_n} w$ then $\rep_J(y) = \rep_J(w)$. \\[0.15em]

(ii) & If $y \sim_{\operatorname{I}_{n-1},\,\operatorname{II}_{n-1} ,\,\operatorname{III}_{n-1}} w$, then $\rep_K(y) = \rep_K(w)$.
\end{tabular}
\end{lemma}


\begin{corollary}\label{corollary 7.6}
Let $y,w \in W_n$. Then

\smallskip

\setlength{\tabcolsep}{4pt}
\hskip -5.5pt \begin{tabular}{rl}

(i) & $w = \ep^k(y)$ for some $k \in \Z$ if and only if $y \sim_{\operatorname{I}_n,\,\operatorname{II}_n} w$, \\[0.15em]

(ii) & $w = \psi^k(y)$ for some $k \in \Z$ if and only if $y \sim_{\operatorname{I}_{n-1},\,\operatorname{II}_{n-1} ,\,\operatorname{III}_{n-1}} w$.
\end{tabular}
\end{corollary}

Proof -- We prove only part (i). Suppose $w = \ep^k(y)$. By the definition of the left extension of a map (see \S\,\ref{S3.1}), we have $\rep_J(y) = \rep_J(w)$. Now set $u:= \pr_J(y)$, $v:=\pr_J(w)$, and $x:=\rep_J(y)= \rep_J(w)$. Then
\[ \begin{array}{cccr}
\exists\, k \in \Z \,:\, v = \ep^k(u) & \iff & u \sim_{R,J} v & \  \text{by property (A4)} \\
& \iff & P(u) = P(v) & \ \text{\S\,5 of \cite{KL}}\\
& \iff & u \sim_{\operatorname{I}_n,\,\operatorname{II}_n} v & \ \text{Theorem~\ref{K Knuth}} \\
& \iff & xu \sim_{\operatorname{I}_n,\,\operatorname{II}_n} xv. & 
\end{array}\]
For the last equivalence, the `left to right' implication is 
given by  Lemma~\ref{type A proj lemma}, while the converse uses Lemma~\ref{lemma rep relations} (i) as well.  \qed

\ms

The following result is due to 
Welsh through private communication.

\begin{lemma}\label{Welsh 2} (Welsh) Let $w =w(1), \ldots,  w(n) \in W_n \setminus \Area_n$ be a signed permutation such that $w(n-1)$ and $w(n)$ have opposing signs. Then the element $w'= w s_{n-1}$ is connected to $w$ by a sequence of Knuth relations of types $\operatorname{I}_n$, $\operatorname{II}_n$ and $\operatorname{III}_{n-1}$.
\end{lemma}



Recall from Corollary~\ref{corollary 7.6} that if $y \sim_{\operatorname{I}_n,\,\operatorname{II}_n} w$, then we have $w = \ep^{k_1} (y)$ for some $k_1 \in \Z$, and if $y \sim_{\operatorname{III}_{n-1}} w$, then we have $w = \psi^{k_2} (y)$ for some $k_2 \in \Z$. So under the hypothesis of Lemma~\ref{Welsh 2}, we know that $w$ can be transformed into $w' = w  s_{n-1}$ by applying a combination of the maps $\ep$ and $\psi$ suitably. In other words, if $w \notin \Area_n$ and $w \sim_{\operatorname{III}_{n}} ws_{n-1}$, then 
$ws_{n-1} \in \orb^R_\Xi(w)$.

\subsection{$\Xi$-orbits} We begin by showing how $W_n \setminus \rArea_n$ is partitioned by the $\Xi$-orbits.
Let $y,w \in W_n$. Then:
\begin{equation}\label{eqn subset}
y  \stackrel\Xi\llra w \ \Ra \ A_n(y) = A_n(w).
\end{equation} 
Indeed, $\Xi$-orbits are subsets of asymptotic right cells by Theorem~\ref{Thm BG 6.2}, and so (\ref{eqn subset}) follows by the characterisation of these cells in Theorem~\ref{ThmBI}~(ii). 

A converse to (\ref{eqn subset}) exists if we assume that $y,w \in W_n\setminus\rArea_n$. Combining our discussion of the elements $w_J$ and $w_Jw_0$ in \S\,\ref{S2.1} with Proposition~\ref{BI Knuth}, it suffices to show that for $y,w \in W_n\setminus\Area_n$, if $y\sim_{\operatorname{I}_{n},\,\operatorname{II}_{n} ,\,\operatorname{III}_{n}} w$ then $y$ and $w$ lie in the same $\Xi$-orbit. However, this follows from Corollary~\ref{corollary 7.6} and Lemma~\ref{Welsh 2}, yielding the following result.

\begin{proposition}\label{orbit iff bitableaux}
Suppose $(W_n,S,\wei)$ is such that $b/a > n-2$, and let $y,w \in W_n \setminus \rArea_n$. Then
\[ y \stackrel\Xi\llra w \iff A_n(y) = A_n(w).\]
\end{proposition}

\ms

We now consider how $\rArea_n$ is partitioned into $\Xi$-orbits. 

\begin{lemma}\label{area orbit iff bitableaux}
Suppose $(W_n,S,\wei)$ is such that $b/a > n-2$, and let $y,w \in \rArea_n$. Then we may have $A_n(y) = A_n(w)$ without 
$w \in \orb_\Xi^R(y)$. 
Indeed, let $1 \les q \les n-1$, and let $T_q$ be any standard bitableau of shape $(1^{n-q} \,|\, 1^q)$. Then the set $\Ga'(T_q):=\{ w \in W_n \, : \, A_n(w) = T_q\}$ is equal to the union of exactly two $\Xi$-orbits.
\end{lemma}

Proof -- 
Consider $w=w(1), \ldots, w(n) \in \Ga'(T_q)$ as a signed permutation, and note that $q=\ell_t(w)$. As in \S\,\ref{S2.3}, we write $x_1, \ldots, x_q$ for the (increasing) subsequence of negative integers in the row form of $w$, and $y_1, \ldots, y_{n-q}$ for the (decreasing) subsequence of positive integers. Theorems 1 and 2 of Schensted \cite{Schensted} motivate the following observations.

\begin{itemize}
\item Suppose $w(n) > 0$; that is, $w(n) = y_{n-q}$. Then by Lemma~\ref{type A proj lemma}, $\pr_J(w)$ has a longest increasing subsequence of length $q+1$, and a longest decreasing subsequence of length $n-q$. By Lemma~\ref{type B proj lemma} (ii), we note that $\ell_t(w) = \ell_t(\pr_K(w))$. \\[-1em]

\item Suppose $w(n)< 0$; that is, $w(n) = x_q$. Then by Lemma~\ref{type A proj lemma}, $\pr_J(w)$ has a longest increasing subsequence of length $q$, and a longest decreasing subsequence of length $n-q+1$. By Lemma~\ref{type B proj lemma} (ii), we note that $\ell_t(w) = \ell_t(\pr_K(w))+1$.
\end{itemize}

We denote $r:=n-q$ and $\ob{x}_i:=|x_i|$ for reasons of typesetting. Applying the generalised Robinson--Schensted insertion algorithm, we see that:
\[\begin{tikzpicture}
\node at (-1.05,0.9) {$A_n(w) \, =$};

\xbox{0}{2.7}{0.9}{y_r}
\xbox{0}{1.8}{0.9}{y_{r-1}}
\zbox{0}{0.9}{0.9}{\vdots}
\xbox{0}{0}{0.9}{y_1}

\xbox{1.5}{2.7}{0.9}{\ob{x}_q}
\xbox{1.5}{1.8}{0.9}{\ob{x}_{q-1}}
\zbox{1.5}{0.9}{0.9}{\vdots}
\xbox{1.5}{0}{0.9}{\ob{x}_1}

\node at (2.65,0.87) {.};
\end{tikzpicture}\]

\begin{figure}
\[\begin{tikzpicture}

\node at (-1.7,0.9)  {$A_{n-1}(\pr_K(w)) \, = $};

\tbox{0}{2.7}{0.9}{y'_{r-1}}
\tbox{0}{1.8}{0.9}{y'_{r-2}}
\zbox{0}{0.9}{0.9}{\vdots}
\tbox{0}{0}{0.9}{y'_{1}}

\tbox{1.5}{2.7}{0.9}{\ob{x}'_q}
\tbox{1.5}{1.8}{0.9}{\ob{x}'_{q-1}}
\zbox{1.5}{0.9}{0.9}{\vdots}
\tbox{1.5}{0}{0.9}{\ob{x}'_1}

\node at (2.7,0.87) {\text{,}};

\node at (6,0.9) {$\begin{array}{lr} \text{where if } w(i)= x_j & \text{ (resp. } y_j \text{),} \\  \text{then }x'_j:=\pr_K(w)(i) & \text{ (resp. } y'_j \text{),} \end{array}$};

\node at (-2.3,-1.65) {$\iff w(n) >0$};

\node at (-1.95,-3) {$ \iff  P(\pr_J(w)) \, = $};

\tbox{0}{-2.2}{0.9}{x''_1}
\tbox{0.9}{-2.2}{0.9}{x''_{2}}
\xbox{1.8}{-2.2}{0.9}{\cdots}
\tbox{2.7}{-2.2}{0.9}{x''_q}
\tbox{3.6}{-2.2}{0.9}{y''_r}
\tbox{0}{-3.1}{0.9}{y''_{r-1}}
\tbox{0}{-4}{0.9}{y''_{r-2}}
\zbox{0}{-4.9}{0.9}{\vdots}
\tbox{0}{-5.8}{0.9}{y''_{1}}

\node at (4.8,-3.2) {\text{,}};

\node at (6,-4.7) {$\begin{array}{lr} \text{where if } w(i)= x_j & \text{ (resp. } y_j \text{),} \\  \text{then }x''_j:=\pr_J(w)(i) & \text{ (resp. } y''_j \text{).} \end{array}$};
\end{tikzpicture}\]

We also have: \qqq \qqq \qqq \qq

\[\begin{tikzpicture}

\node at (-1.7,0.9)  {$A_{n-1}(\pr_K(w)) \, = $};

\tbox{0}{2.7}{0.9}{y'_r}
\tbox{0}{1.8}{0.9}{y'_{r-1}}
\zbox{0}{0.9}{0.9}{\vdots}
\tbox{0}{0}{0.9}{y'_{1}}

\tbox{1.5}{2.7}{0.9}{\ob{x}'_{q-1}}
\tbox{1.5}{1.8}{0.9}{\ob{x}'_{q-2}}
\zbox{1.5}{0.9}{0.9}{\vdots}
\xbox{1.5}{0}{0.9}{\ob{x}'_{1}}

\node at (2.7,0.87) {\text{,}};

\node at (6,0.9) {$\begin{array}{lr} \text{where if } w(i)= x_j & \text{ (resp. } y_j \text{),} \\  \text{then }x'_j:=\pr_K(w)(i) & \text{ (resp. } y'_j \text{),} \end{array}$};

\node at (-2.25,-1.6) {$\iff w(n) < 0$};

\node at (-1.9,-3) {$ \iff P(\pr_J(w)) \, = $};

\tbox{0}{-2.15}{0.9}{x''_1}
\tbox{0.9}{-2.15}{0.9}{x''_{2}}
\xbox{1.8}{-2.15}{0.9}{\cdots}
\tbox{2.7}{-2.15}{0.9}{x''_q}
\tbox{0}{-3.05}{0.9}{y''_{r}}
\tbox{0}{-3.95}{0.9}{y''_{r-1}}
\zbox{0}{-4.85}{0.9}{\vdots}
\tbox{0}{-5.75}{0.9}{y''_{1}}

\node at (3.9,-3.15) {\text{,}};

\node at (6,-4.65) {$\begin{array}{lr} \text{where if } w(i)= x_j & \text{ (resp. } y_j \text{),} \\  \text{then }x''_j:=\pr_J(w)(i)& \text{ (resp. } y''_j \text{).} \end{array}$};
\end{tikzpicture}\]
\caption{\small Two collections of equivalent statements pertaining to the proof of Lemma~\ref{area orbit iff bitableaux}.}\label{Lemma 7.10 diagram}
\end{figure}

\noi Now we use our observations to deduce the contents of Figure~\ref{Lemma 7.10 diagram}.

\ms

We have a disjoint union of non-empty sets $\Ga'(T_q) \eq \ga'_1 \sqcup \ga'_2$, where 
\begin{itemize}
\item  $\ga'_1:= \{\, w \in W_n \, : \, A_n(w) = T_q \text{\ \qs and\  \qs } w(n)>0 \,\}$, \\[-1em] 
\item $\ga'_2:=\{\, w \in W_n \,:\, A_n(w) = T_q \text{\ \qs and\ \qs } w(n) <0 \,\}$. 
\end{itemize}
By noting that $\Ga'(T_q)$ coincides with an asymptotic right cell, we see that this decomposition may be identified with that in (a right-handed version of) Corollary~\ref{subcells}. 
Let $y,w \in \Ga'(T_q)$. Then using Figure~\ref{Lemma 7.10 diagram}, we deduce that the following four statements are true if and only if  $y,w \in \ga'_i$ for some $i \in \{1,2\}$:
\[ \begin{array}{rllrl}
(i) & A_{n-1}(\pr_K(y))=A_{n-1}(\pr_K(w)), &  &(ii) & y \stackrel{\{(K,\psi)\}}\llra w, \\
(iii) & P(\pr_J(y)) = P(\pr_J(w)), & & (iv) & y \stackrel{\{(J,\ep)\}}\llra w.
\end{array}\]

\noi This proves the lemma, as well as Proposition~\ref{asy orbits} (iii). \qed


\begin{proposition}\label{asy orbits}
Suppose $(W_n,S,\wei)$ is such that $b/a > n-2$. \\[-0.8em]

\setlength{\tabcolsep}{4pt}
\hskip -5.5pt \begin{tabular}{r p{318pt}}
(i) & $W_n$ comprises $\operatorname{YBT} +\, 2^n -2$ \qs $\Xi$-orbits, where $\operatorname{YBT}$ is the number of standard bitableaux of size $n$. \\[0.5em]

(ii) & Let $w \in W_n$.
\end{tabular}
\begin{itemize}[leftmargin=36pt]
\item If $w \in W_n \setminus \rArea_n$ then
\[ \orb^R_\Xi(w) \eq \{\,w' \in W_n \,:\, A_n(w) = A_n(w')\,\}. \qq \qq \qq \ \]
\item If $w \in \rArea_n$ and $w(n)>0$, then
\[ \orb^R_\Xi(w) \eq  \{\,w' \in W_n \,:\, A_n(w) = A_n(w') \text{ and } w'(n)>0\,\}.\]
\item If $w \in \rArea_n$ and $w(n)<0$, then
\[ \orb^R_\Xi(w) \eq \{\,w' \in W_n \,:\, A_n(w) = A_n(w') \text{ and } w'(n)<0\,\}.\]
\end{itemize}

\setlength{\tabcolsep}{4pt}
\hskip -8.5pt \begin{tabular}{r p{318pt}}
(iii) & If $w \in \Area_n$, then 
\end{tabular}
\[ \orb^R_\Xi(w) \eq \orb^R_{\{(J,\ep)\}}(w) \eq  \orb^R_{\{(K,\psi)\}}(w).\]

\end{proposition}
Proof -- Part (ii) follows from Proposition~\ref{orbit iff bitableaux} and the proof of Lemma~\ref{area orbit iff bitableaux}, 
while part (i) holds by part (ii) and application of (\ref{area number of cells 2}) from \S\,\ref{S2.1}.  \qed

\ms

An immediate consequence of this is that $w \in W_n\setminus\rArea_n$ if and only if the $\Xi$-orbit containing $w$ coincides with the asymptotic right cell containing~$w$.

\ms

As noted in \S\,\ref{3.notation}, $\Xi$ satisfies the conditions in Lemma~\ref{KL pair lemma}; together with this result we confirm the conjecture mentioned at the end of Example~6.8 of \cite{BG}.

\ms

A table comparing the number of $\Xi$-orbits to the number of asymptotic and intermediate right cells of $W_n$ is given in Figure~\ref{orbit table}.

\begin{figure}
\[\begin{array} {|c|cc|}
 \hline
 \phantom{\Big{|}}  & \frac ba = n-1  & \frac ba > n-1  \\ [0.2 em] \hline
\phantom{\Big{|}}\!\!n\!=\!2\!\!\phantom{\Big{|}} & 4 & 6  \\
\phantom{\Big{|}}3\phantom{\Big{|}} & 16 &  20   \\
\phantom{\Big{|}}4\phantom{\Big{|}} & 68 & 76   \\
\phantom{\Big{|}}5\phantom{\Big{|}} & 296 &  312  \\
\phantom{\Big{|}}6\phantom{\Big{|}} & 1352 &  1384  \\
\phantom{\Big{|}}7\phantom{\Big{|}} & 6448 & 6512 \\[0.2 em] \hline
\end{array} \ \
\begin{array}{|c|}
\hline
\frac ba > n-2  \phantom{\Big{|}}   \\ \hline
\phantom{\Big{|}}8\phantom{\Big{|}} \\
\phantom{\Big{|}}26\phantom{\Big{|}}\\
\phantom{\Big{|}}90\phantom{\Big{|}}  \\
\phantom{\Big{|}}342\phantom{\Big{|}} \\
\phantom{\Big{|}}1446\phantom{\Big{|}} \\
\phantom{\Big{|}}6638\phantom{\Big{|}} \\ \hline
\end{array}\]

\caption{On the left, the number of left/right cells of $W_n$ with respect to the given weights. On the right, the number of $\Xi$-orbits that comprise $W_n$.}\label{orbit table}
\end{figure}

\subsection{The property ($\star$)} 
We require a final discussion regarding bipartitions and $\Xi$-orbits before moving on to Vogan $(\Xi,\rxi)$-classes.

\ms

If $\pi$ is a partition of $n$ with conjugate denoted $\pi'$, and $\lam = (\lam^+ \,|\, \lam^-) = (\lam^+_1, \ldots, \lam^+_f \,|\, \lam^-_1, \ldots, \lam^-_g) \Vdash n$ is a bipartition of $n$, then $\lam' := ( (\lam^-)' \,|\, (\lam^+)' )$ is the conjugate bipartition of $\lam$.
If we write $\lam = (\lam_1, \ldots, \lam_{f+g})$ where $\lam_i=\lam^+_i$ for $1 \les i \les f$ and $\lam_{f+j} = \lam^-_j$ for $1 \les j \les g$, then we can define:
\[ I(\lam):= S \setminus \{t,\, s_{\lam_1},\, s_{\lam_1+\lam_2},\,  \ldots,\,   s_{\lam_1 + \lam_2 + \cdots + \lam_{f+g-1}}  \}. \]
Set $q:= \lam_1 + \cdots + \lam_f$, denote by $w_{I(\lam)}$ the longest word of $W_{I(\lam)}$, and by $w_q$ the longest word of $W_q$. Then we set
\[ w_\lam \, := \,  w_{I(\lam)} w_q.\]
Thus $\rdes(w_{I(\lam)}) = I(\lam)$, and if $b/a > q-1$, then $\rxi(w_\lam) \cap \{t_1, \ldots, t_n\} = \{t_1, \ldots, t_q\}$. With the notation of \S\,\ref{S5}, we have $w_{(q \,|\, n-q)} = \si_{n,q} \in \Omega_{\ze_q}$. 

Note that $\sh(w_\lam)=\lam'$, and so $w_{\lam'} \in \Omega_\lam$. We are interested in whether or not an element $z \in W_n$ has the property:
\[(\star) \qq \qq \qq \ \  \ \ \ \ \exists\, \nu \in \mathcal{V}_\Xi \ : \ \nu(z) \in \orb^L_\Xi (w_{\sh(z)'}).  \qq \qq \qq \qq \  \  \ \ \]
For such an element $z$, it follows from Proposition~\ref{asy orbits} (ii) that there exists some $\hat{z} \in \orb^L_\Xi(z)$ such that $\nu(\hat{z}) = w_{\sh(z)'}$; this $\hat{z}$ is the unique element in the intersection of $\orb^L_\Xi(z)$ with $\orb^R_\Xi(w_{\sh(z)'})$. 

\ms

\begin{remark}\label{orbit scope}
From Proposition~\ref{asy orbits}, we can state that $w \in W_n$ satisfies ($\star$) if and only if $w \in W_n\setminus\rArea_n$ or $w \in \rArea_n$, $w(n)>0$ and $w^{-1}(n) > 0$.
\end{remark}

\begin{remark}\label{rmk tsuki}
Let $y,w \in W_n\setminus \rArea_n$. By Propositions~\ref{orbit iff bitableaux} and~\ref{BI Knuth} respectively, we have:
\[ y \stackrel\Xi\llra w \iff A_n(y) = A_n(w) \iff y \sim_{\operatorname{I}_n,\,\operatorname{II}_n ,\,\operatorname{III}_n}w.\]
\end{remark}


\section{Vogan $(\Xi,\rxi)$-classes when $b/a > n-2$}

In this section, we prove the main result of the paper. To do this, we begin by stating two results of Bonnaf\'{e} from \cite{BonKLB}. 
 Recall that $\chi_q := s_{n-1} \cdots s_q$.  

\begin{proposition}\label{Bon Prop 4.1} (Bonnaf\'{e}, Proposition 4.1 of \cite{BonKLB})
Let $0 \les q \les n$, let $\al, \be \les_e p_{n,q}$ and let $\si, \si' \in W_{S\setminus\{t,s_q\}}$ be such that $\si \sim_L \si'$. Assume that $b/a \ges n-1$. Then
\[ \al a_q \si \be^{-1} \,\sim_L\, a_q \si' \be^{-1}.\]
\end{proposition}

\begin{proposition}\label{Bon Prop 6.1} (Bonnaf\'{e}, Proposition 6.1 of \cite{BonKLB})
Let $0 \les q \les n-1$ and assume that $b/a \les n-1$. Then
\[ s_1 \cdots s_q \cdot \si_{n,q} \,\sim_L\, \si_{n,q+1} \chi_{q+1}^{-1}.\]
\end{proposition}

Note that these two results of Bonnaf\'{e}'s are valid simultaneously if and only if $b/a = n-1$. 

\begin{lemma}\label{Lemma 001}
Suppose $(W_n,S,\wei)$ is such that $b/a \ges n-1$ and $y,w \in \Area_n$. Then the following are equivalent: \\[-0.7em]

\setlength{\tabcolsep}{4pt}
\hskip -5.5pt \begin{tabular}{rl}
(i) & $y \app^{\Xi,\rxi} \!\!\! w$, \\[0.1em]
(ii) &$y \sim_L w$, \\[0.2em]
(iii) & $\rxi(y)=\rxi(w)$.
\end{tabular}
\end{lemma}

Proof -- We have (ii) implies (i) by Theorem~\ref{Thm BG 7.2}, while Definition~\ref{def vogan} shows that (i) implies (iii). If $b/a > n-1$, then (iii) implies (ii) by Lemma~\ref{lemma 5.10}~(i). If $b/a = n-1$, then we have $\rxi(y)=\rxi(w)$ if and only if $\rdes(y)=\rdes(w)$ by Corollary~\ref{corollary 5.9} (ii). So, it remains to prove that when $b/a=n-1$,  $\rdes(y)=\rdes(w)$ implies that $y \sim_L w$.

\ms

We have seen in Lemma~\ref{lemma 5.10} (iii) and the subsequent discussion in \S\,\ref{S2.6} that 
$\Up(w):=\{z \in \Area_n \,:\, \rdes(z)=\rdes(w)\}$ 
comprises exactly two asymptotic left cells, a decomposition which we denote by $\Up(w) = \Ga_1 \sqcup \Ga_2$. Our first step will be to show that
\begin{equation}\label{eqn 09}
y \sim_L w \text{ \qs when \qs } \frac ba > n-1 \ \qs \Ra \qs \ y \sim_L w  \text{ \qs when \qs } \frac ba = n-1.
\end{equation}
Following our discussion of the elements $w_J$ and $w_Jw_0$ in \S\,\ref{S2.1}, we can assume $y, w \in \rArea_n$ when proving (\ref{eqn 09}).

Suppose $y,w \in \Ga \subseteq \rArea$, where $\Ga$ is an asymptotic left cell. We can decompose $\Ga$ as $\Ga = \ga_1 \sqcup \ga_2$, where the disjoint union is as in Corollary~\ref{subcells}. Let $q=\ell_t(y)=\ell_t(w)$. By Corollary~\ref{subcells}, the sets $\ga_1$ and $\ga_2$ are subsets of left cells with respect to any choice of parameters; in particular, they must be subsets of left cells when $b/a = n-1$. Thus to prove (\ref{eqn 09}), it suffices to show that there exist representatives $y' \in \ga_1$ and $w' \in \ga_2$ such that $y' \sim_L w'$ when $b/a = n-1$. So set $y':=\si_{n,q} \theta^{-1}$ where $\theta \les_e p_{n,q}$ is such that $y' \in \Ga$, and set $w':= \chi_q y' \in \ga_2$. Let $b_q$ be as in \S\,\ref{S2.4}.

 We now apply Proposition~\ref{Bon Prop 4.1} with $\al = \chi_q$, $\be = \theta$ and $\si=\si'=b_q$ to see that $y' \sim_L w'$ when $b/a = n-1$; (\ref{eqn 09}) is now proved.

\ms

By Proposition~\ref{intermediate area cells} (i), we know that the sets $\Ga_1$ and $\Ga_2$ are of the form $\Ga(\si_{n,q}\tau^{-1})$ and $\Ga(\si_{n,q+1}\chi^{-1}_{q+1}\tau^{-1})$ respectively for some $0 \les q \les n-1$ and $\tau \les_e p_{n-1,q}$, which we fix for the remainder of this proof.  Thus it remains to show that $\si_{n,q}\tau^{-1} \sim_L \si_{n,q+1}\chi^{-1}_{q+1}\tau^{-1}$ when $b/a = n-1$. Indeed, the following holds for $b/a \in (n-2,n-1]$, resulting in Corollary~\ref{proto-areadiff}.


By Corollary~\ref{subcells}, we have $\si_{n,q} \sim_L s_1\cdots s_q \si_{n,q}$ for all choices of parameters. Then by Proposition~\ref{Bon Prop 6.1}, we have $\si_{n,q} \sim_L \si_{n,q+1}\chi_{q+1}^{-1}$. Recall that if $y \sim_L w$, then $\nu(y) \sim_L \nu(w)$ for all $\nu \in \mathcal{V}_\Xi$. 

Consider $\Ga'(\si_{n,q})$, and the corresponding subset $\ga_1$ as in (a right-handed version of) Corollary~\ref{subcells}. Compare $\ga_1$ with the description of $\Xi$-orbits in Proposition~\ref{asy orbits} (ii) to see that
\[ \si_{n,q} \stackrel\Xi\llra \si_{n,q} \tau^{-1} \text{ if and only if } \tau \les_e p_{n-1,q}.\]
Now we use  Proposition~\ref{asy orbits} (iii) to reduce this statement to 
\[ \si_{n,q} \stackrel{\{(I,\de)\}}\llra \si_{n,q} \tau^{-1} \text{ if and only if } \tau \les_e p_{n-1,q},\]
where $(I,\de) \in \Xi$. By the same results, we have 
\[ \si_{n,q+1}\chi_{q+1}^{-1}  \stackrel{\{(I,\de)\}}\llra \si_{n,q+1}\chi_{q+1}^{-1} \tau^{-1} \text{ if and only if } \tau \les_e p_{n-1,q}.\]
It may be verified that $\pr_I(\si_{n,q}) = \pr_I(\si_{n,q+1}\chi_{q+1}^{-1})$ for $(I,\de) \in \Xi$. It now follows that $\nu ( \si_{n,q}) = \si_{n,q} \tau^{-1}$ if and only if $\nu(\si_{n,q+1}\chi_{q+1}^{-1}) = \si_{n,q+1}\chi_{q+1}^{-1} \tau^{-1}$. Thus $\si_{n,q}\tau^{-1} \sim_L \si_{n,q+1}\chi^{-1}_{q+1}\tau^{-1}$, as needed. \qed

\begin{corollary}\label{proto-areadiff}
No asymptotic left cell contained in $\Area_n$ is a left cell with respect to $b/a \in (n-2,n-1]$.
\end{corollary}

We will write $\app$ for $\app^{\Xi,\rxi}$ in proofs for cleaner notation. Remark~\ref{orbit scope} describes the scope of the following result.

\begin{lemma}\label{Lemma 002a}
Suppose $(W_n,S,\wei)$ is such that $b/a > n-2$. Let $y,w \in W_n$ be elements satisfying 
($\star$) with $\ell_t(y) = \ell_t(w)$. Then
\[ y \app^{\Xi,\rxi} \!\! w \ \Ra \ y \in \orb^L_\Xi(w).\]
\end{lemma}

We note that the following proof is adapted from part of Ariki's proof of Theorem A in \S\,3.4 of \cite{Ariki}, who in turn remarks to have adapted an argument from the proof of Satz 5.25 in Jantzen \cite{Jantzen}. 

\ms

Proof -- Denote $q:= \ell_t(y) = \ell_t(w)$. Let $\lam,\mu \Vdash n$ be such that $\sh(y) = \lam'$ and $\sh(w) = \mu'$. Let $\al,\be \in \mathcal{V}_\Xi$ be such that $\al(y) \in \orb^L_\Xi(w_\lam)$ and $\be(w) \in \orb^L_\Xi(w_\mu)$. Let $\hat{y}, \hat{w} \in W_n$ be the unique elements in the intersection of  $\orb^L_\Xi(y)$ with $\orb^R_\Xi(w_\lam)$ and $\orb^L_\Xi(w)$ with $\orb^R_\Xi(w_\mu)$ respectively. 
Note that $\al(\hat{y}) = w_\lam$ and $\be(\hat{w}) = w_\mu$. 

As $\hat{y} \in \orb^L_\Xi(y)$ and $\hat{w} \in \orb^L_\Xi(w)$, we have $y \sim_L \hat{y}$ and $w \sim_L \hat{w}$. Applying Theorem~\ref{Thm BG 7.2}, we see that 
\[\hat{y} \,\app\, y \,\app\, w \,\app\, \hat{w}.\]

We apply the generalised Robinson--Schensted algorithm to the first $q$ entries of the row forms of $\al(\hat{w})$ and $\be(\hat{y})$; that is, their negative entries. Our focus will first be on the recording bitableaux of $\al(\hat{w})$, with that of $\be(\hat{y})$ being obtained in a similar fashion.

As $\hat{y} \app \hat{w}$ and $\al(\hat{y}) = w_\lam$, we have $\rxi(\al(\hat{w})) = \rxi(\al(\hat{y})) = \rxi(w_\lam)$. 
So by Lemma~\ref{signdes}, the first $\lam^+_1$ entries in the row form of $ \al(\hat{w})$ form an increasing sequence (of negative integers), the next $\lam^+_2$ entries form an increasing sequence, and so on. 
So the first column of the second tableau of $B_n(\al(\hat{w}))$ is a column of length at least $\lam^+_1$. As $\sh(\al(\hat{w})) = \sh(w) = \mu'$, the length of this first column equals $\mu^+_1$. 

Similarly, the length of the first column of the second tableau of $B_n(\be(\hat{y}))$ is a column of length at least $\mu^+_1$, and equal to $\lam^+_1$. Combining these observations indicates that $\lam^+_1 = \mu^+_1$. 

We continue in this way, noting that bumping between the columns during the algorithm produces a contradiction to the previous statements. We thus determine that $\lam^+ = \mu^+$. By similarly working with the remaining $n-q$ entries in the row forms of $\al(\hat{w})$ and $\be(\hat{y})$ -- that is, the positive entries -- we determine that $\lam^- = \mu^-$. Thus $\lam = \mu$, and we have:
\[ B_n(\al(\hat{w})) \eq B_n(w_\lam)  \eq  B_n(w_\mu) \eq B_n(\be(\hat{y})).\]

As $\be(\hat{y}) \in \orb^R_\Xi(w_\lam)$, we have $A_n(\be(\hat{y})) = A_n(w_\lam)$ by Proposition~\ref{asy orbits}~(ii). Similarly, $\al(\hat{w}) \in \orb^R_\Xi(w_\mu)$, and so $A_n(\al(\hat{w})) = A_n(w_\mu)$. Now that we know the recording and insertion bitableaux of $\al(\hat{w})$ and $\be(\hat{y})$, we can determine these elements; we have:
\[ \al(\hat{w}) \eq \be(\hat{y}) \eq w_\lam \eq w_\mu \eq \al(\hat{y}) \eq \be(\hat{w}).\]
Thus $\hat{y} = \hat{w}$, and so $y \in \orb^L_\Xi(w)$ as needed. \qed

\begin{lemma}\label{Lemma 002b}
Suppose $(W_n,S,\wei)$ is such that $b/a > n-2$.  Let $y,w \in W_n\setminus\Area_n$ with $\ell_t(y) \neq \ell_t(w)$. Then $y$ and $w$ do not lie in the same left Vogan $(\Xi,\rxi)$-class. 
\end{lemma}

Proof -- Suppose, towards a contradiction, that $y \app w$. As $\rxi(y) = \rxi(w)$, we must have $\ell_t(y),\qs \ell_t(w) \in \{n-1,n\}$ and $b/a \in (n-2,n-1]$. We may assume that $\ell_t(y) = n-1$ and $\ell_t(w) = n$ without loss of generality.

As $y$ satisifies ($\star$) by Remark~\ref{orbit scope}, there exists some $\al \in \mathcal{V}_\Xi$ such that $\al(y) \in \orb^L_\Xi(w_{\sh(y)'})$. Thus $\rxi(\al(y)) \cap \{t_1, \ldots, t_n\} = \{t_1, \ldots, t_{n-1}\}$. As $\ell_t(y) = n-1$, we have $\al(y)(n-1) < 0 < \al(y) (n)$ by Lemma~\ref{signdes}. 
By Remark~\ref{rmk tsuki}, 
$\al(y)$ and $\al(y)\cdot s_{n-1}$ lie in the same $\Xi$-orbit. 
 So let $\ga \in \mathcal{V}_\Xi$ be such that 
$\ga(y) = \al(y) s_{n-1}$.  We now have $\ga(y)(n-1) > 0 > \ga(y)(n)$, and so 
\[\rxi(\ga(y)) \cap \{t_1, \ldots, t_n\} \eq \{t_1, \ldots, t_{n-2}\}.\] 
On the other hand, as $\ell_t(w)=n$, we have 
\[ \rxi(\nu(w)) \cap \{t_1, \ldots, t_n\} \eq \{t_1, \ldots, t_{n-1}\} \ \text{ for all } \nu \in \mathcal{V}_\Xi.\] 
As $\rxi(\ga(y)) \neq \rxi(\ga(w))$ implies that $y \napp w$, a contradiction is reached. \qed

Finally, we need to show that $\Area_n$ is closed under $\app^{\Xi,\rxi}$\!\!.

\begin{lemma}\label{Lemma 003}
Suppose $(W_n,S,\wei)$ is such that $b/a > n-2$. Let $y \in \Area_n$ and $w \in W_n$. Then
\[y \app^{\Xi,\rxi} \!\! w \ \Ra \ w \in \Area_n.\]
\end{lemma}

Proof -- If $w$ does not satisfy ($\star$), then $w \in \rArea_n$ by Remark~\ref{orbit scope}, and we are done. So we may assume from now on that $w$ satisfies ($\star$). 

Suppose that $y$ satisfies ($\star$). By Lemma~\ref{Lemma 002a}, we may assume that $\ell_t(y) \neq \ell_t(w)$;  thus $\ell_t(y),\,\ell_t(w) \in \{n-1,n\}$ and $b/a \in (n-2,n-1]$. 

If $\ell_t(y) = n-1$, then $y = \si_{n,n-1}$. If $\ell_t(y) = n$, then $y = \si_{n,n}$. In either case, we have 
$\rxi(w) = \rxi(y) = \{t_1,\ldots, t_{n-1}\}$.
Application of Lemma~\ref{signdes} then implies that $w \in \Area_n$. 
So we may assume from now on that $y$ does not satisfy ($\star$). By Remark~\ref{orbit scope}, we have $y \in \rArea_n$, and thus $\ell_t(y) \les n-1$.

So let $\lam,\mu \Vdash n$ be such that $\sh(y) = \lam'$ and $\sh(w) = \mu'$, and let $\be \in \mathcal{V}_\Xi$ be such that $\be(w) \in \orb^L_\Xi(w_\mu)$. This implies that $\rxi(\be(w)) = \rxi(w_\mu)$. 

\ms

Case A -- Suppose that $\ell_t(y) = \ell_t(w)$; denote this value by $q$. 

As $y\app w$ and $\be(w)(n) > 0$, it must be the case that $\be(y)(n)>0$. By Proposition~\ref{asy orbits} (ii), we have:
\begin{equation}\label{y(n)>0}
\nu(y)(n) > 0 \qq \ \text{for all } \,\nu \in \mathcal{V}_\Xi.
\end{equation}
Suppose now, towards a contradiction, that $w\notin \Area_n$. By Remark~\ref{rmk tsuki}, we may apply a particular sequence of generalised Knuth relations of type $\operatorname{III}_n$ to $\be(w)$ while remaining in $\orb^R_\Xi(w)$; we have: 
\[\be(w)\  \stackrel\Xi\llra\  \be(w) \cdot s_q  \ \stackrel\Xi\llra \ \ \cdots \ \ \stackrel\Xi\llra \ \be(w) \cdot s_q s_{q+1} \cdots s_{n-1}.\]

\noi So let $\ga \in \mathcal{V}_\Xi$ be such that $\ga(w) = \be(w) \cdot s_q s_{q+1} \cdots s_{n-1}$. Then we have: 
\begin{itemize}
\item $\be(w) ( j) < 0 \iff j \in \{1, \ldots, q\}$  \qs and \\[-1.1em]
\item $\ga(w) ( j) < 0 \iff j \in \{1, \ldots, q-1,n\}$. \\[-1.1em]
\end{itemize}
As $y \app  w$ and $q \les n-1$, we also have $\ga(y)(j)< 0$ if and only if $j \in \{1, \ldots, q-1, n\}$. However, this contradicts (\ref{y(n)>0}).

\ms

Case B -- Suppose that $\ell_t(y) \neq \ell_t(w)$. As $\rxi(y) = \rxi(w)$ and $\ell_t(y) \les n-1$, we have $\ell_t(y) = n-1$, $\ell_t(w) = n$ and $b/a \in (n-2,n-1]$. Note that
\begin{itemize}
\item $\rxi(y) \cap \{t_1, \ldots, t_{n}\} = \rxi(w) \cap \{t_1, \ldots, t_{n}\} = \{t_1, \ldots, t_{n-1}\}$,  \\[-1.0em]
\item $\rxi(\si_{n,n-1}) \cap \{t_1, \ldots, t_{n}\} = \{t_1, \ldots, t_{n-1}\}$.
\end{itemize} 
As $y,\qs \si_{n,n-1} \in \Area_n$, by Corollary~\ref{corollary 5.9} (i) we have $\rxi(w) = \rxi(\si_{n,n-1})$. We know that $\ell_t(w) = n$, and so we have enough information to apply Lemma~\ref{signdes} and determine that $w = \si_{n,n} \in \Area_n$. \qed 


\begin{theorem}\label{Thm2}
Suppose $(W_n,S,\wei)$ is such that $b/a \ges n-1$. Then
\[y \sim_L w \iff y \app^{\Xi,\rxi} \!\! w.\]
\end{theorem}

Proof -- This follows from Lemmas \ref{Lemma 001}, \ref{Lemma 002a}, \ref{Lemma 002b} and \ref{Lemma 003}. \qed

\ms


Using the results of this section, we may now give a combinatorial description of the left Vogan $(\Xi,\rxi)$-classes for $b/a > n-2$. The classes differ depending on whether $b/a > n-1$ or $b/a \in (n-2,n-1]$, and they do so precisely on $\Area_n$.

If $b/a > n-1$, then we obtain a description of $\Area_n$ using Proposition~\ref{asy cells in area}. If $b/a \in (n-2,n-1]$,  it is sufficient to determine the classes for a representative in this interval. So take $b/a = n-1$, and look to the decomposition of $\Area_n$ given in Proposition~\ref{intermediate area cells} (ii). 

\newpage

\begin{corollary}\label{hmm, corollary}
If $y,w \in W_n\setminus\Area_n$ and $b/a > n-2$, then
\[y^{-1} \stackrel\Xi\llra w^{-1} \iff y \sim_L w \iff B_n(y)=B_n(w) \iff y \app^{\Xi,\rxi} \!\!w.\]
If $b/a \in (n-2,n-1]$ we have:
\[\Area_n \eq\bigsqcup_{q{\scalebox{0.6}{$\,$}}={\scalebox{0.6}{$\,$}}0}^{n-1}\ \bigsqcup_{\tau\, \les_e\, p_{n-1,q}} \sm \{\, \pi \cdot\si_{n,q}\cdot \tau^{-1} \,:\, \pi \les_e \Pi_{n,q}=s_{n-q-1} \cdots s_1 \cdot t \cdot p_{n,q}\,\},\]
with each set being a left Vogan $(\Xi,\rxi)$-class, and a left cell if $b/a = n-1$.

\ms

\noi If $b/a > n-1$ we have:
\[\Area_n \eq\bigsqcup_{q{\scalebox{0.6}{$\,$}}={\scalebox{0.6}{$\,$}}0}^{n}\ \bigsqcup_{\tau\, \les_e\, p_{n,q}}  \! \{\, \pi \cdot\si_{n,q}\cdot \tau^{-1} \,:\, \pi \les_e  p_{n,q}\,\},\]
with each set being a left Vogan $(\Xi,\rxi)$-class and a left cell.

\ms

\noi Further, for $b/a > n-2$, $\Area_n$ is closed under the relation $\app^{\Xi,\rxi}$\!\!.
\end{corollary}

From this, we may infer some information about the intermediate two-sided cells of $W_n$.

\begin{corollary}
Let $\Omega \subseteq W_n \setminus \Area_n$ be an asymptotic two-sided cell. Then $\Omega$ is contained within an intermediate two-sided cell. Further, $\Area_n$ is contained within an intermediate two-sided cell of $W_n$.
\end{corollary}

We may also characterise the asymptotic left cells that change when passing to the intermediate or sub-asymptotic case.

\begin{corollary}\label{areadiff}
Let $\Ga \subseteq W_n$ be an asymptotic left cell. Then $\Ga$ is also a left cell with respect to $b/a \in (n-2,n-1]$ if and only if $\Ga \subseteq W_n\setminus\Area_n$.
\end{corollary}


\subsection*{Acknowledgements}
 This paper has its origins in the author's thesis. As such, it would not exist without the patient supervision of both Lacri Iancu and Jean-Baptiste Gramain. Additional thanks are due to Meinolf Geck for his invaluable comments. 
\bibliography{mybib}{}
\bibliographystyle{plain}

\end{document}